\documentclass[authoryear,times,3p]{elsarticle}

\usepackage{caption}

\usepackage[labelformat=simple,belowskip=0em,aboveskip=0pt]{subcaption}    

\usepackage{makecell}

\usepackage{fancyvrb}

\usepackage{amssymb}
\usepackage{amsthm}

\theoremstyle{definition}
\newtheorem{theorem}{Theorem}[section]

\newtheorem{example}[theorem]{Example}

\usepackage{siunitx}

\usepackage{hyperref}

\usepackage{float}
\newfloat{mylistt}{tbph}{lol}[section]
\floatname{mylistt}{List}

\usepackage{newfloat}
\DeclareFloatingEnvironment[
fileext = lol,
placement = htbp,
name = List]{mylist}

\usepackage{enumitem}

\usepackage{CJKutf8}

\usepackage{amsmath}
\DeclareMathOperator*{\argmax}{argmax}
\DeclareMathOperator*{\argmin}{argmin}

\usepackage{todonotes}
\usepackage{mathtools}

\usepackage{bbm}

\usepackage{placeins}

\newcommand{\dcos}{d_\mathrm{{cos}}}
\newcommand{\norm}[1]{\lVert{#1}\rVert}

\newcommand{\datadep}[1]{{{#1}}}

\newcommand{\bitstream}[2]{\mathrm{Bitstream}(#1,#2)}

\DeclareMathOperator{\im}{im}

\journal{Somewhere}

\begin{document}

\begin{frontmatter}



\title{A topological analysis of the space of recipes}


\author[inst1]{Emerson G. Escolar}
\ead{e.g.escolar@people.kobe-u.ac.jp}
\ead[url]{https://emerson-escolar.github.io}

\author[inst1]{Yuta Shimada}
\ead{246d502d@stu.kobe-u.ac.jp}

\author[inst1]{Masahiro Yuasa}
\ead{yuasa@people.kobe-u.ac.jp}

\affiliation[inst1]{
  organization={Kobe University Graduate School of Human Development and Environment},
  addressline={3-11 Tsurukabuto, Nada},
  city={Kobe City},
  postcode={657-8501},
  state={Hyogo},
  country={Japan}}

\begin{abstract}
  In recent years,
  the use of data-driven methods has provided insights into
  underlying patterns and principles behind culinary recipes.
  In this exploratory work, we introduce the use of topological data analysis, especially persistent homology,
  in order to study the space of culinary recipes.
  In particular, persistent homology analysis provides a set of recipes surrounding
  the multiscale ``holes'' in the space of existing recipes.
  We then propose a method to generate novel ingredient combinations using
  combinatorial optimization on this topological information.
  We made biscuits using the novel ingredient combinations,
  which were confirmed to be acceptable enough
  by a sensory evaluation study.
  Our findings indicate that topological data analysis has
  the potential for providing new tools
  and insights in the study of culinary recipes.
\end{abstract}




\end{frontmatter}



\section{Introduction}
\label{sec:sample1}

The application of data science to food is a novel and rapidly growing field
with the potential to answer many interesting questions in gastronomy
\citep{ahnert2013network,mouritsen2017data,bagler2018data,min2019survey,herrera2021contribution,goel2022computational}.
In this study, we introduce the application of new ideas from topological data analysis
in the study of cooking recipes.
In particular, persistent homology \citep{landi1997new,frosini1999size,robins1999towards,edelsbrunner2002topological},
one of the main tools in topological data analysis (see for example the books \citep{edelsbrunner2010computational,carlsson2021topological,dey2022computational}),
is able to describe the multi-scale
connected components, holes, voids (cavities), and so on, in data.
The hope is that by understanding the shape of the data via topological quantities describing it,
one can get hints as to the data generating process behind
the phenomenon \citep{chazal2021introduction}.
An initial impetus for this work is to apply this novel data analysis method
(and related methods) in order to further understand the ``shape'' of the 
space of cooking recipes.
Furthermore, once ``holes'' can be identified in the recipe space,
these may point to potential opportunities for creating novel recipes.

\paragraph{Contributions}

In this exploratory work, we introduce the use of topological data analysis,
especially persistent homology, in order to study the space of culinary recipes.
For this, we provide both an intuitive description with examples in the main text,
and also provide mathematically precise definitions in the Appendix.
Persistent homology analysis provides a set of recipes surrounding
the multiscale ``holes'' in the space of existing recipes,
which we exploit using combinatorial optimization in order to
generate novel ingredient combinations.
Our method optimizes the dissimilarity of the ingredient combination
compared to existing recipes,
while at the same time maintaining coherence by using the topological information
consisting of a cycle of recipes detected by persistent homology.
We perform analysis showing that the ingredient combinations generated by our method
are novel with respect to the input data.
We also selected some related ingredient combinations obtained from the analysis
and confirmed that we can cook variations of a dish (cream cheese biscuits) from it.
A sensory evaluation study confirmed that these biscuits were acceptable enough.
We also discuss some limitations of our study and the potential for further explorations
in Section~\ref{sec:discussion}.

\subsection{Related literature}

  \paragraph{Network science and topological ideas}
Our use of topological data analysis is related to
prior work applying network analysis to the study of food recipes \citep{herrera2021contribution}.
One groundbreaking work in this field is the paper \citep{ahn2011flavor},
approaching the problem of identifying general principles of
``food-pairing'' \citep{blumenthal2008big,briscione2018flavor,coucquyt2020art}.
That is, they studied whether or not there are any general principles that support
the ingredient combinations used in existing recipes,
versus the astronomical number of theoretical possible ingredient combinations.
The work \citep{ahn2011flavor} approaches the topic using mathematical techniques from
complex networks, to encode relations between ingredients and flavors as a network,
called the \emph{flavor network}
(nodes are ingredients and edges are determined by shared flavor compounts between ingredients).
Via their analysis, they argue that Western cuisines tends to prefer recipes
that contains ingredients with similar flavor compounds,
while Eastern Asian cuisines tend to avoid flavor sharing \citep{ahn2011flavor}.
Further investigations into food-pairing have also been performed to other regional cusines
\citep{varshney2013flavor,jain2015analysis,issa2018analysis,tallab2016exploring,al2021exploring,dougan2023computational},
and an extension called food bridging \citep{simas2017food} has been proposed.

Following \citep{ahn2011flavor}, we also focus our attention on the food recipes,
but for this initial proof-of-concept study, we do not consider the flavor profiles of different ingredients,
but instead focus on combinations of ingredients.
We believe that our approach gives a complementary view on the recipe space,
but a further study involving flavor profiles would also be interesting.

A similar approach has been adopted by \citep{kular2011using}.
In that study, they constructed a network of $300$ recipes from many different cuisines,
where the nodes are recipes and edges are determined by the shared ingredients between recipes.
They found that the network exhibits small-world and scale free properties.
Furthermore, for some non-European cusines, they \citep{kular2011using} were able to detect
those cuisines using network community analysis.

\paragraph{Persistent homology}

In the field of text mining and image analysis,
the paper \citep{wagner2014towards}
uses persistent homology to reveal a global structure of similarities 
in their context of image and text documents.
Furthermore, they argue that it is essential to incorporate higher-dimensional relationships,
and suggests that topological data analysis can ``point to parts of the feature-space which are not populated.''
This is the perspective we adopt,
as we use the topological methods to find
regions (``holes'') of the recipe space we could exploit to suggest novel recipes.
In particular, we take persistent homology as a tool for exploratory data analysis.

Another line of application of persistent homology
is in materials science \citep{hiraoka2016hierarchical,saadatfar2017pore,kimura2018non}.
There, the interpretation of ``holes''
as detected by persistent homology is straightforward,
as the data is usually a set of points (or an image) in $2$D or $3$D space,
and a hole corresponds
directly with a ``physical'' hole that one can visualize as-is\footnote{Note that (persistent) homology describes a hole
  by an equivalence class of cycles surrounding it, and thus the representative cycles that can be used
  to illustrate it are not unique. This issue is discussed further in Subsection~\ref{subsec:ph},
  and is independent of the point being raised here,
  which simply says that once a representative cycle for a hole is chosen and fixed,
  it can be directly drawn in $2$D or $3$D.}.
In contrast,
our target space of the recipe space 
is high-dimensional,
and thus similar visualizations are not as straightforward.

\paragraph{Machine learning and recipes}
In recent years, the application of machine learning techniques for
research into various tasks related to recipes and food.
As a detailed survey and complete review of literature is beyond the scope of this paper,
we refer the reader to survey papers
\citep{min2019survey,min2019food,goel2022computational,bondevik2023systematic}.
%
We especially note a few papers that are more closely related to our
theme of ingredient combinations.
For example, \citep{kazama2018neural}
proposes a system that can quantify the ``regional cuisine style'' mixture
of recipes and can suggest ingredient substitutions to transform a recipe
to fit better to a selected regional cuisine style.
The work \citep{park2019kitchenette} uses Siamese neural networks to propose
a model called ``KitcheNette'' that predicts food ingredient pairing scores
and can be used to recommend ingredient pairings.
Furthermore,
\citep{park2021flavorgraph}
incorporates large-scale information about flavor molecules in food
ingredients as a ``FlavorGraph'' and uses that information to
recommend ingredient pairings.

\paragraph{Recipes and chef’s creativity}

Recipes are important tools for any food and culinary culture,
especially for preparing dishes in a delicious, nutritious, and safe manner \citep{borghini2015recipe}.
In general, recipes consist mainly of
food ingredients and their quantities or weights,
cooking methods and time,
kitchen utensils and equipment, and
food presentation.
In general, recipes are employed for various reasons and in various situations,
for example
when cooking a dish for the first time,
cooking traditional cuisines,
for cooking practice or culinary training,
for regulating nutritional qualities of dishes\footnote{For example, it is indicated that online diet management apps are useful for health management in type 2 diabetes patients \citep{tominaga2022individuals,hironaka2024impact}.},
for reducing cooking time and cost,
and so on.
In addition, various types of media exist for recipes,
such as  books, magazines, web pages, and apps.
%
%
In this work, due to limitations in data and to simplify our analysis,
we take a simplified view of recipes as combinations of ingredients;
see the subsection~\ref{subsec:data} on data.

In general, chefs are experts who can use many recipes and cook many
different dishes.
%
Especially, creative chefs are regarded as culinary artists \citep{ekincek2023recipe},
because a chef’s culinary creativity is defined as the artistic and novel expression of the chef’s inner world,
and the creativity is transferred to cuisine \citep{lee2021pate,lee2020creative}.
Ekincek and Günay (2023) 
suggested that creative chefs are influenced by internal factors such as
memories, instincts, travels, and experiences
as well as external factors such as
family, physical environment, culture, professional environment, and art \citep{ekincek2023recipe}.
However,
the variety of dishes that a chef can create
depends on the chef’s creativity, which can be subjective
and potentially limited by the range of internal and external factors.
A tool for suggesting novel ideas for dishes may be beneficial
for augmenting a creative chef's potential for novelty,
and the use of tools from data science may contribute to this.


\section{Materials and Methods}
\label{sec:matmet}

\subsection{Data}
\label{subsec:data}

We obtained recipe data from Supplementary Dataset 2 of \citep{ahn2011flavor}.
The data\footnote{
Furthermore, this is different from the ingredient -- flavor compound data
used in the same work \citep{ahn2011flavor} from which they constructed the flavor network
(nodes are ingredients and edges are determined by shared flavor compounts between ingredients).}
is a list of recipes, where
each recipe is given as a list (set) of ingredients used, and a
``region''\footnote{For this dataset, the regions are `African', `EastAsian', `EasternEuropean', `LatinAmerican',
`MiddleEastern', `NorthAmerican', `NorthernEuropean', `SouthAsian',
`SoutheastAsian', `SouthernEuropean', and `WesternEuropean'} associated to that recipe.
The data contains duplicated ingredient lists, possibly caused by differently-named recipes using the same ingredients; these duplicates are removed.
The basic statistics of this cleaned data are described in
Table~\ref{table:datastats}.
\begin{table}[h]
  \centering
  \begin{tabular}{r|r}
    Number of recipes & \datadep{48,983} \\
    Number of ingreds.\ & \datadep{381} \\
    Ave.~num.~ingreds.~per recipe & \datadep{8.4936}\\
    Std.~dev.~num.~ingreds.~per recipe & \datadep{3.5091}
  \end{tabular}
  \caption{Basic Statistics}
  \label{table:datastats}
\end{table}

In the network science perspective (see for example \citep{herrera2021contribution}),
this data forms a \emph{recipe--ingredient bipartite network},
where the nodes are recipes and ingredients
and edges are formed by linking each recipe to the ingredients it uses.
Focusing on relationships between ingredients,
one can consider the network with ingredients as nodes
and with edges determined by whether or not two ingredients appear together in some recipe.
This forms the \emph{ingredient co-occurence network}.
Instead, our focus is on the recipes themselves, and we consider recipes as the nodes in our network.
We discuss this construction below.

Each recipe is considered as a $0$-$1$ vector via one-hot encoding of ingredients.
In particular, we first order the list of all possible (relevant) ingredients $i_1, i_2, \hdots, i_M$.
Then, a recipe is considered as the vector $x \in \mathbb{R}^M$,
where the $j$th coordinate $x_j$ is $1$ if the recipe uses
ingredient $i_j$, and $0$ otherwise.
Throughout, we shall freely identify (sub)sets with $0$-$1$ vectors.
We let $X$ be the set of all the $0$-$1$ vectors of the recipes.

We use the cosine similarity (as used in the field of information retrieval; see for example \citep{singhal2001modern}),
to measure the similarity of any two given recipes. In general, for vectors $x, y \in \mathbb{R}^M$, the cosine similarity between $x$ and $y$ is defined to be
\[
  s_\text{cos}(x,y) := \frac{x \cdot y}{\norm{x}\norm{y}} = \mathrm{cos}(\theta)
\]
where $\theta$ is the angle between $x$ and $y$.
For use with persistent homology, instead of similarity, we need a dissimilarity measure,
so we use the \emph{cosine dissimilarity} defined as
\[
  \dcos(x,y) := 1 - s_\text{cos}(x,y) =
  1 - \frac{x \cdot y}{\norm{x}\norm{y}}.
\]
In general, the more the vectors $x$ and $y$ are pointing in the same direction,
$\dcos$ becomes closer to $0$,
while if $x$ and $y$ are closer to pointing in opposite directions,
$\dcos$ becomes closer to its maximum possible value of $2$.
For our data, since recipes $x$ and $y$ have non-negative entries,
the maximum possible value of $\dcos(x,y)$ is $1$,
attained when $x$ and $y$ are orthogonal ($x \cdot y = 0$).
We also remark that, given the way we are encoding recipes as vectors, the inner product
$x \cdot y$ is exactly equal to the number of ingredients shared between recipes $x$ and $y$.
%
In Appendix~\ref{subsec:appendix:dissim} we perform further analysis
of the distribution of dissimilarities in the data.
\begin{figure}[h]  
  \centering
  \includegraphics[width=0.45\textwidth]{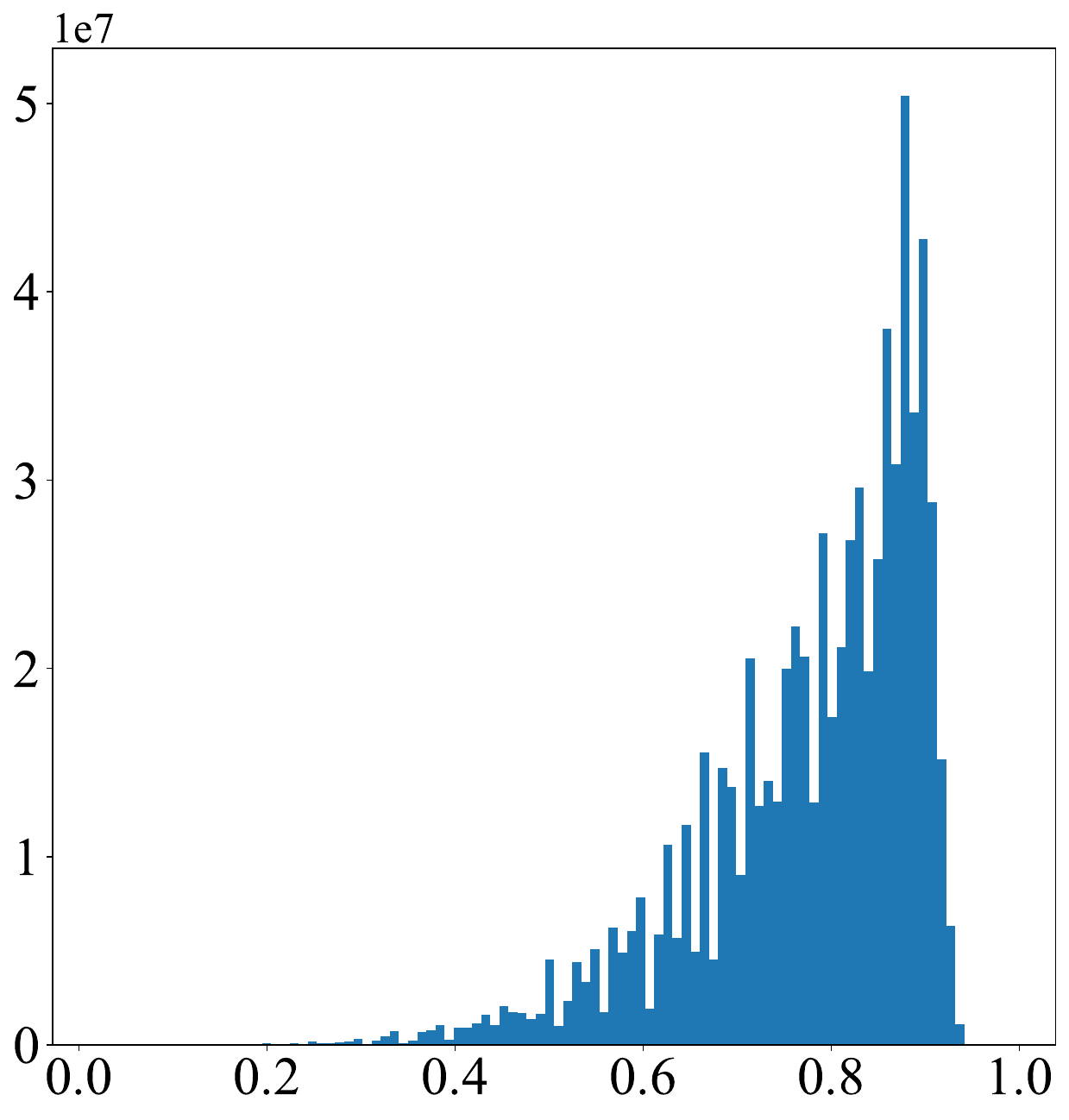}
  \caption{Histogram of cosine dissimilarities 
    of distinct pairs of distinct recipes from the recipe data \citep{ahn2011flavor},
    restricted to dissimilarities less than $1$.
    There were \datadep{482,978,610} (out of \datadep{1,199,642,653}) pairs
    with dissimilarity exactly equal to $1$. These are the pairs of recipes
    that share no ingredients at all.
  }
  \label{fig:histodists}
\end{figure}

\subsection{The Vietoris-Rips complex}
\label{subsec:vr}

We review basic ideas from topological data analysis.
See for example the books \citep{edelsbrunner2010computational,carlsson2021topological,dey2022computational}.

We associate a geometric shape to the recipe data  using the
mathematical construction known as the Vietoris-Rips complex,
which is an example of a simplicial complex.
The idea is the following. Choose some threshold value $t$.
For each pair of recipes, if they are sufficiently similar
(that is, if their dissimilarity is at most the chosen threshold)
then we connect them by an edge. Similarly, for every triple of recipes,
if each pair in the triple is sufficiently similar then we create a triangle, and so on.

In general,
an \emph{abstract simplicial complex} on a set of nodes (also called \emph{vertices})
$X = \{x_1,x_2, \hdots, x_n\}$ is a set $K$ of nonempty subsets of $X$
satisfying the condition that $\sigma \in K$ and $\emptyset \neq \tau \subset \sigma$ implies that $\tau \in K$.
The elements $\sigma \in K$ (i.e.\ $\sigma \subset X$) are called the \emph{simplices} of $K$.
The simplices of $K$ can be interpreted in the following way.
A singleton $\{v\}$ in $K$ is just the node $v$,
a two-element set $\{x,y\}$ in $K$ can be thought of as an edge connecting $x$ and $y$,
a three-element set $\{x,y,z\}$ in $K$ can be thought of as a triangle involving $x$, $y$, $z$,
and so on.
The \emph{dimension} of a simplex $\sigma$ is its number of elements minus $1$.
A vertex is $0$-dimensional, an edge $1$-dimensional, and so on.
A simplex with dimension $q$ will be also called a $q$-simplex.
Each nonempty $\tau \subset \sigma$ is called a \emph{face} of the simplex $\sigma$.
For example, the faces of triangle $\sigma = \{x,y,z\}$ are:
$\sigma$ itself, the three edges $\{x,y\}$, $\{y,z\}$, $\{x,z\}$, and the three vertices
$\{x\}$, $\{y\}$, $\{z\}$.
Furthermore, one can think of a simplicial complex as a higher-dimensional generalization of a graph.
In fact, a simple graph (with no loops and multiple edges) with vertex set $X$ and edge set $E$
can be thought of as a simplicial complex (with $K = \{ \{v\} \mid v \in X \} \cup E$).

For a set $X$, a \emph{dissimilarity} on $X$ is a
function $d: X \times X \rightarrow \mathbb{R}_{\geq 0}$ satisfying $d(x,x) = 0$ and $d(x,y) = d(y,x)$
for all $x,y \in X$.
A pair $(X,d)$ of a finite set $X$ together with a dissimilarity $d$ on $X$
is called a (finite) \emph{dissimilarity space}.
Then, for a dissimilarity space $(X, d)$,
its \emph{Vietoris-Rips complex}\footnote{
  We note that the usual definition of a Vietoris-Rips complex
  assumes that $d$ is a metric; that is, 
  $d(y,x) = d(x,y)$ for all $x, y \in X$, with $d(x,y) = 0$ if and only if $x=y$,  
  and the triangle inequality ($d(x,z) \leq d(x,y) + d(y,z)$ for all $x,y,z \in X$) is satisfied.
  However, as also noted in \citep{chazal2014persistence}, the definition continues to make sense for
  a dissimilarity $d$, as we have presented here.
  We also remark that the cosine dissimilarity is not a metric, because the triangle inequality may not
  hold for certain choices of $x,y,z$.
}
at scale $t$, denoted $V_t(X)$, is the
abstract simplicial complex
containing all nonempty subsets $\sigma$ of $X$ such
that every pair $x,y$ in $\sigma$ has dissimilarity at most $t$,
or in symbols:
\[
  V_t(X) :=
  \left\{
    \emptyset \neq \sigma \subset X  \mid d(x,y) \leq t \text{ for all } x,y \in \sigma
  \right\}.
\]
This definition captures the above intuition of forming edges, triangles, etc.\
based on comparing pairwise dissimilarities against the threshold value $t$.

\begin{example}
  \label{example:vrcomplex}
  Let us consider a small example of a Vietoris-Rips complex $K = V_t(X)$ for the set
  $X = \{x_1,x_2,x_3,x_4,x_5\} \subset \mathbb{R}^4$ with
  \begin{align*}
    x_1 &= [1, 0, 0, 1] \\
    x_2 &= [0, 0, 1, 1] \\
    x_3 &= [0, 1, 1, 0] \\
    x_4 &= [1, 1, 0, 0] \\
    x_5 &= [1, 0, 0, 0]
  \end{align*}
  For example, this could arise from a data set of five recipes $x_1,\hdots,x_5$
  involving four ingredients $i_1,\hdots,i_4$,
  say, $i_1 = \text{coffee}$, $i_2 = \text{milk}$, $i_3 = \text{sugar}$, $i_4 = \text{cinammon powder}$
  where $x_1$ corresponds to the ``recipe'' using coffee and cinammon powder,
  $x_2$ the ``recipe'' of sugar and cinammon power, and so on.

  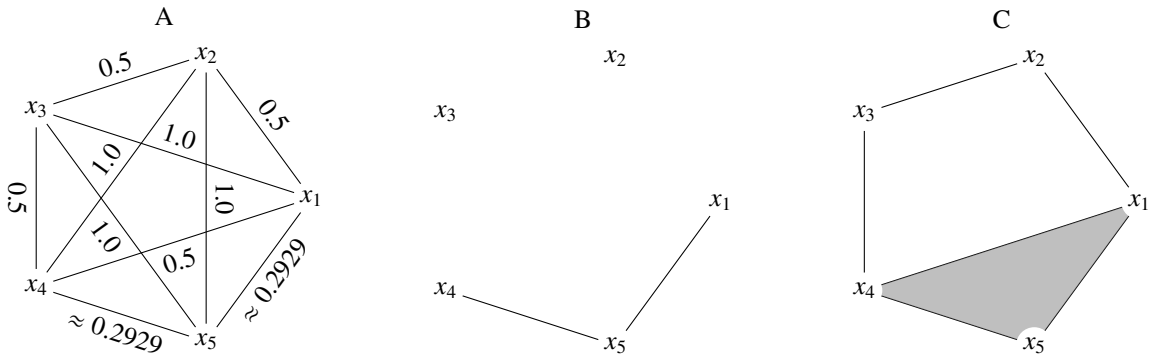
\begin{figure}[h!]
    \captionsetup[sub]{font=normalsize}
    \centering
    \begin{subfigure}[b]{0.33\textwidth}    
      \centering
      \caption{}
      \label{subfig:exampleVRdists}
      \begin{tikzpicture}[
        line join = round, line cap = round, baseline=(current bounding box.center)]
        \tikzset{main node/.style={circle,fill=black, inner sep=1pt},label node/.style={scale=0.75},}  
        \coordinate (1) at (0:2);
        \coordinate (2) at (72:2);
        \coordinate (3) at (144:2);
        \coordinate (4) at (216:2);
        \coordinate (5) at (288:2);

        \draw (1) edge node[above,sloped] {0.5} (2);
        \draw (1) edge node[above,sloped] {1.0} (3);
        \draw (1) edge node[below,sloped] {0.5} (4);
        \draw (1) edge node[below,sloped] {$\approx 0.2929$} (5);
        \draw (2) edge node[above,sloped] {0.5} (3);
        \draw (2) edge node[above,sloped] {1.0} (4);
        \draw (2) edge node[above,sloped] {1.0} (5);
        \draw (3) edge node[below,sloped] {0.5} (4);
        \draw (3) edge node[below,sloped] {1.0} (5);
        \draw (4) edge node[below,sloped] {$\approx 0.2929$} (5);
        
        \foreach \x in {1,2,3,4,5} {
          \draw (\x) node [main node,fill=white]{$x_\x$};
        }
      \end{tikzpicture}
    \end{subfigure}
    \hfill
    \begin{subfigure}[b]{0.33\textwidth}    
      \centering
      \caption{}
      \label{subfig:exampleVRt1}
      \begin{tikzpicture}[
        line join = round, line cap = round, baseline=(current bounding box.center)]
        \tikzset{main node/.style={circle,fill=black, inner sep=1pt},label node/.style={scale=0.75},}  
        \coordinate (1) at (0:2);
        \coordinate (2) at (72:2);
        \coordinate (3) at (144:2);
        \coordinate (4) at (216:2);
        \coordinate (5) at (288:2);  

        \draw (4) -- (5) -- (1);
        

        \foreach \x in {1,2,3,4,5} {
          \draw (\x) node [main node,fill=white]{$x_\x$};
        }
      \end{tikzpicture}
    \end{subfigure}
    \hfill
    \begin{subfigure}[b]{0.33\textwidth}    
      \centering
      \caption{}
      \label{subfig:exampleVRt2}
      \begin{tikzpicture}[
        line join = round, line cap = round, baseline=(current bounding box.center)]
        \tikzset{main node/.style={circle,fill=black, inner sep=1pt},label node/.style={scale=0.75},}  
        \coordinate (1) at (0:2);
        \coordinate (2) at (72:2);
        \coordinate (3) at (144:2);
        \coordinate (4) at (216:2);
        \coordinate (5) at (288:2);  

        \draw (4) -- (5) -- (1) -- (2) -- (3) -- cycle;
        \draw (1) -- (4);
        \begin{scope}[]
          \draw[draw=none,fill=gray,fill opacity=.5] (1)--(4)--(5)--cycle;    
        \end{scope}
        
        \foreach \x in {1,2,3,4,5} {
          \draw (\x) node [main node,fill=white]{$x_\x$};
        }
      \end{tikzpicture}    
    \end{subfigure}
    \caption{
      Illustrations for the Vietoris-Rips complex in Example~\ref{example:vrcomplex}.
      \subref{subfig:exampleVRdists},
      cosine dissimilarities between pairs of points.
      \subref{subfig:exampleVRt1},
      Vietoris-Rips complex $V_{t_1}(X)$ with threshold $t_1 = 1-\frac{1}{\sqrt{2}}$.      
      \subref{subfig:exampleVRt2},
      Vietoris-Rips complex $V_{t_2}(X)$ with threshold $t_2 = 0.5$.}      
    \label{fig:exampleVR}
  \end{figure}

  In Figure~\ref{subfig:exampleVRdists}, 
  we show the pairwise cosine dissimilarities between points of $X$ as labels on edges,
  where for example
  \[
    \dcos(x_4,x_5) = \dcos(x_1,x_5) = 1 - \frac{1}{\sqrt{2}} \approx 0.2929.
  \]
  Then, for example
  $V_0(X)$ (with threshold $t=0$) is simply the set of points $X$,
  and
  $V_{t_1}(X)$ with threshold $t_1 = 1-\frac{1}{\sqrt{2}}$ is
  as illustrated in Figure~\ref{subfig:exampleVRt1} and
  $V_{t_2}(X)$ with threshold $t_2 = 0.5$
  is as illustrated in Figure~\ref{subfig:exampleVRt2}.
\end{example}


We compare and contrast this with a common way of constructing a network from recipes
(for example the ``Network of Recipes'' construction of \citep{kular2011using}, see also \citep{herrera2021contribution}),
where nodes are recipes and each pair of recipes $r_1$ and $r_2$
is given an edge with weighting equal to the number of shared ingredients between $r_1$ and $r_2$
(we can delete the edges with weighting $0$, i.e.\ no shared ingredients).
Note that when recipes are encoded as $0$-$1$ vectors with respect to the ingredients,
the number of shared ingredients can also be expressed by the inner product $r_1 \cdot r_2$.
In contrast, our construction uses the cosine dissimilarity
$
\dcos(r_1,r_2)  =
  1 - \frac{r_1 \cdot r_2}{\norm{r_1}\norm{r_2}}.
$
which normalizes the inner product and turns it into a dissimilarity measure.
Furthermore, the Vietoris-Rips complex considers higher-order relations
(triples of recipes, quadruples, and so on) instead of just pairs of recipes.

\subsection{(Persistent) homology}
\label{subsec:ph}

Here, we provide a quick introduction to the basic ideas and definitions for homology and persistent homology \citep{edelsbrunner2002topological}.
For any $q \in \{0,1,2,\hdots\}$, we construct below a vector space called the
``$q$th-degree homology group with coefficients in $\mathbb{F}_2$'', denoted $H_q(K)$,
of a simplicial complex $K$ (for example, the Vietoris-Rips complex $K = V_t(X)$ at a fixed scale $t$).
We give precise definitions in Appendix~\ref{sec:defns}.

Intuitively, $H_q(K)$ represents the $q$-dimensional ``holes'' (features) in $K$,
where its dimension as a vector space counts the number of independent features. For example,
$H_0(K)$ counts the number of connected components of $K$,
$H_1(K)$ the number of independent loops,
$H_2(K)$ the number of independent cavities,
and so on.
In Example~\ref{example:vrcomplex},
the $1$st (i.e.~$q=1$) homology group of the Vietoris-Rips complex $K = V_{t_2}(X)$
with threshold $t_2 = 0.5$ is the vector space generated by the homology class of the cycle
\[
  z_1 = \{x_3,x_2\} + \{x_2,x_1\} + \{x_1,x_4\} + \{x_4,x_3\}
\]
which represents the simple cycle that goes around the vertices $x_3,x_2,x_1,x_4$ and back to $x_3$.
Such a $z_1$ is called a \emph{representative cycle}, and the set of cycles \emph{equivalent}
to it (i.e.\ cycles that differ from it by the boundary of a set of triangles) is called its homology class $[z_1]$.
This homology class $[z_1]$ can also be given by the representative
\[
  z_2 = \{x_3,x_2\} + \{x_2,x_1\} + \{x_1,x_5\} + \{x_5,x_4\} + \{x_4,x_3\}
\]
which takes a detour through $x_5$. Their difference is the boundary of the triangle $\{x_1,x_5,x_4\}$,
and thus $[z_1] = [z_2]$; at the homological level, these two cycles are considered the same.

For the Vietoris complex, what is the correct threshold $t$ to use?
Persistent homology bypasses this question by considering all possible (relevant) choices of $t \in \mathbb{R}$
simultaneously, as follows.
First, we note that for $t \leq t'$, $V_t(X) \subseteq V_{t'}(X)$.
This means that $V_t(X)$ forms a \emph{filtration} of simplicial complexes (an increasing sequence of
simplicial complexes).
Persistent homology tracks the appearance (birth) and disappearance (death)
of homology classes as we increase $t$.

A more precise statement is given in Appendix~\ref{sec:defns},
and the conditions below refer to the conditions in Theorem~\ref{theorem:pd}. For each $q \in \{0,1,2,\hdots\}$,
persistent homology of the Vietoris-Rips filtration $\{V_t(X)\}_{t}$
provides a multiset of pairs of thresholds $\{(b_i,d_i)\}_{i=1}^s$
together with a set of $q$-cycles $\{c_i\}_{i=1}^s$, satisfying the following conditions.
First, for each $i \in \{1,2,\hdots,s\}$,
the cycle $c_i$ is created exactly at threshold $t = b_i$ (conditions \ref{ph:b1} and \ref{ph:b2}).
Then, looking at its homology class at threshold $t$,
$[c_i]_t$ is nonzero for $t$ satisfying $b_i \leq t < d_i$ (condition \ref{ph:alive}),
and becomes zero at threshold $t = d_i$ and beyond (condition \ref{ph:d}).
Finally, the choice of the set $\{c_i\}$ and pairs $\{(b_i,d_i)\}_{i=1}^s$
should capture all independent features at every scale $t$ (condition \ref{ph:basis}).

The multiset of pairs is called the
$q$th-degree persistence diagram of $\{V_t(X)\}_{t}$, which we shall denote by
$D_q(X) = \{(b_i,d_i)\}_{i=1}^s$.
While the multiset $\{(b_i,d_i)\}_{i=1}^s$ is uniquely defined,
the set of cycles $\{c_i\}_{i=1}^s$ is not.
We call $c_i$ a \emph{representative cycle}
for the \emph{birth-death pair} $(b_i,d_i)$.
For a birth-death pair $(b_i,d_i)$, its \emph{lifespan} is the value $d_i - b_i$.

A more direct interpretation can be given in our chosen setting
(see also the intepretation in \citep[Section~7]{wagner2014towards}),
especially for $q=1$.
Consider a $q$-cycle $c$ associated to the birth-death pair $(b,d)$.
A $q$-cycle is composed of $q$-simplices; with $q=1$, this is just a collection of edges (pairs of recipes).
At the birth scale $t = b$ (when the cycle $c$ is just created),
each pair of recipes in the cycle have dissimilarity at most $b$, as a threshold.
As we increase the threshold, more and more dissimilar recipes (equivalently, less and less similar recipes)
are considered. Eventually, at threshold $d$, the death scale,
enough
$2$-simplices (i.e.~triplets of recipes) with dissimilarity at most $d$ are created,
filling-in the cycle.

Furthermore, while we do not use it directly in our work,
it is worth mentioning the following stability result
for persistence diagrams of the Vietoris-Rips filtration of dissimilarity spaces.
It guarantees that persistence diagrams are stable under small perturbations of the dissimilarities.
\begin{theorem}[{\citep[Theorem~5.2~and~Subsection~4.2.5]{chazal2014persistence}, also stated in \citep{turner2019rips}}]
  Let $(X,d_X)$ and $(Y,d_Y)$ be finite dissimilarity spaces. Then,
  \[
    d_I(D_q(X), D_q(Y)) \leq 2 d_{\text{GH}}(X,Y)
  \]
  where $d_{\text{GH}}$ is the Gromov-Hausdorff distance (see \citep[Definition~7.3.10]{burago2001course}).
\end{theorem}

\subsection{Recombination and simplification of identified recipes}
\label{subsec:optimization}

Let us consider the $1$st-degree ($q=1$) persistence diagram and the representative cycles
of the Vietoris complex complex of the recipe data $(X,\dcos)$.
Each representative cycle $c$ is a $1$-cycle, and thus is a
sum of $1$-simplices (edges) $e_1 + e_2 + \hdots + e_p$ with boundary equal to $0$.
The set of vertices (existing recipes) $\{r_1,r_2,\hdots,r_\ell\}$
incident to these edges are the important recipes for this representative cycle
detected by persistent homology.

In the case that $c$ is a simple cycle, 
the cycle $c$ is in the form $c = \{r_1,r_2\} + \{r_2,r_3\} + \hdots + \{r_{\ell-1}, r_\ell\} + \{r_\ell,r_1\}$,
forming a ring around a hole (technically, at least one hole) in $V_b(X)$ at the birth threshold $t=b$.
In the case that $c$ is not a simple cycle
(see Figure~\ref{fig:nonsimple} in Appendix~\ref{subsec:appendix:cycles}), we can still consider the recipes
$\{r_1,r_2,\hdots,r_\ell\}$ involved in the representative cycle $c$ as-is for the analysis below,
or decompose $c$ into simple cycles and use those instead.
For our analysis we perform the latter procedure.

Next, we want to generate a new recipe
$r$ that is located somewhere in the hole surrounded by the cycle $c$.
Naively, we can consider the centroid, in the following way.
Each recipe $r_i$ is a $0$-$1$ vector in $\mathbb{R}^M$. The centroid of the points of $c$ is therefore
\[
  \bar{c} := \frac{1}{\ell}\sum_{i=1}^\ell r_i.
\]
However, $\bar{c}$ is not a $0$-$1$ vector, and its set of nonzero entries $S$ corresponds simply to
the union of all the ingredients in the recipes $r_i$.
Furthermore, the size of $S$ tends to be larger\footnote{Average of \datadep{$36.6$} for the top $10$ cycles with the longest lifespans (see Table~\ref{table:bdl}) versus the average of \datadep{$8.49$} for the dataset.}
than the number of ingredients used in actual recipes,
and so it is difficult to actually create a dish using all the ingredients in $S$.

Instead, we consider the following combinatorial optimization problem.
Given an integer $\nu \geq 2$, we seek a size-$\nu$ subset $y$ of $S$ that is as far (dissimilar)
as possible from existing recipes $X$.
We want a solution
\begin{equation}
  \label{equation:mastercombiopti}
  y_\ast := \argmax_{y \subset S,~|y| = \nu} \dcos(y, X),
\end{equation}
where the dissimilarity is
$
\dcos(y,X) := \min_{x \in X} \dcos(y,x).
$
Note that we identify subsets $y \subset S \subset \{1,\hdots, M\}$ with $0$-$1$ vectors in $\mathbb{R}^M$.
This is equivalent to finding solutions $y_\ast$ for the following min-max problem:
\begin{equation}
  \label{eq:minmax}
  \argmin_{y \subset S,~|y| = \nu} \max_{x \in X} \frac{y \cdot x}{\norm{y}\norm{x}}
\end{equation}
(each $y_\ast$ is a recipe whose \emph{similarity}
with existing recipes is as small as possible).

To reduce the dimensionality, we restrict vectors to $S$, so that we can consider $0$-$1$ vectors in $\mathbb{R}^{|S|}$ instead of in $\mathbb{R}^M$. Without loss of generality, let $S = \{1,2,\hdots, s\}$.
For a vector $x = [x_1,x_2,\hdots,x_M]^T \in \mathbb{R}^M$, let $\pi_S(x) := [x_1,x_2,\hdots,x_{s}]^T \in \mathbb{R}^{s}$.
Then, for $y \subset S$,
(a $0$-$1$ vector in $\mathbb{R}^M$ whose support is contained in $S$),
$\norm{\pi_S(y)} = \norm{y}$,
and
$\pi_S(y) \cdot \pi_S(x) = y \cdot x$ for all $x \in \mathbb{R}^M$.
Furthermore, for  $0$-$1$ vectors $y$ with exactly $\nu$ entries equal to $1$,
$\norm{y} = \sqrt{\nu}$. Thus,
\[
  \frac{y \cdot x}{\norm{y}\norm{x}} =
  \pi_S(y) \cdot  \frac{\pi_S(x)}{\sqrt{\nu}\norm{x}}.
\]
For cases with the number of possible combinations $\binom{|S|}{\nu}$ is small enough,
and with the above dimension reduction,
a brute-force search to solve for $y_\ast$ in Problem~\ref{eq:minmax} may be enough.

We can also further transform Problem~\ref{eq:minmax}. Using the epigraph trick,
this problem is equivalent to solving the following optimization problem
\begin{equation}
  \label{eq:epigraph}
  \argmin\left\{\lambda \;\middle|\; y \subset S,~|y| = \nu, \lambda \in \mathbb{R}, \lambda \geq \frac{x \cdot y}{\norm{x}\norm{y}} \text{ for all } x \in X \right\}.
\end{equation}
Applying the dimensionality reduction to Problem~\eqref{eq:epigraph},
we obtain the following mixed-integer linear programming problem,
where we let 
$X_{S,t} := \left\{\frac{\pi_S(x)}{\sqrt{\nu}\norm{x}} \;\middle|\; x \in X\right\} \subset \mathbb{R}^s$
and
where $\mathbbm{1}$ is the vector with all entries equal to $1$.
\[
  \begin{array}{lll}
    \text{minimize}  & \lambda &\\
    \text{subject to}& v^T y - \lambda \leq 0, & \text{for all }v \in X_{S,t} \\
                     & \mathbbm{1}^T y = \nu, \\
                     & y \in \{0,1\}^s, \\
                     & \lambda \in \mathbb{R}.
  \end{array}
\]
Note that there can be multiple optimal solutions,
representing ties in the dissimilarity to existing recipes.
After obtaining one optimal solution $z$, we add the condition
$
z^T y \leq \nu - 1 
$
to the linear programming problem, which excludes $z$
from the set of feasible solutions without adding
any unnecessary constraints.
Then, we can solve the problem again
and get a different optimal solution.
We repeat this process as necessary.
We use GLPK (GNU Linear Programming Kit, Version 5.0,
\url{http://www.gnu.org/software/glpk/glpk.html})
to solve the mixed-integer linear programming problems.

In essence, the above method:
\begin{enumerate}
\item uses the existing recipes $\{r_1,r_2,\hdots,r_\ell\}$ involved in a 
  representative cycle $c$ as identified by the persistent homology analysis as a starting point,
\item considers the set $S$ of ingredients used in these existing recipes, 
\item and sets up an optimization problem
  to find a (hopefully) new combination of $\nu$ ingredients in $S$ that is as dissimilar as possible to existing recipes.
\end{enumerate}

Thus, the method optimizes the dissimilarity of the ingredient combination
compared to existing recipes,
while at the same time maintaining coherence by using the topological information
consisting of the cycle of recipes detected by persistent homology.
We call a combination of ingredients obtained from our method
a ``solution'' or a ``suggestion'' or a ``suggested combination''.

\subsection{Biscuits preparation}
\label{subsec:biscuitsprep}

To confirm the potential viability of the suggested combinations from our method,
four solutions (from the representative cycle with the longest lifespan)
were selected.
The four suggested combinations consisted of
whole grain wheat flower, starch (corn starch), cranberry (dried cranberry), raisin, gin, and cream cheese (see List~\ref{list:somesolutionsfortopcycle}),
suggesting biscuits.
The composition of the biscuits prepared is presented in Table~\ref{table:biscuitcompo}.
Sugar was added because the recipe data \citep{ahn2011flavor}
does not contain seasonings (e.g. salt and sugar) as ingredients.
Water was added to the No gin (NG) biscuit to equalize the dough weight with other biscuits.
Corn starch was selected as starch because it is used
commonly
in baked goods.    
The control biscuit was prepared using all the ingredients. 
\begin{table}[h]
  \centering  
  \begin{tabular}{l|ccccc}
    \hline
    \thead[b]{Ingredients (g)}
    & \thead[b]{Control\\ { }}
    & \thead[b]{No corn starch\\ (NCS)}
    & \thead[b]{No raisin\\(NR)}
    & \thead[b]{No gin\\ (NG)}
    & \thead[b]{No cranberry\\ (NCB)}\\
    \hline
    Sugar & 30 & 30 & 30 & 30 & 30 \\
    Whole grain wheat flour & 45 & 90 & 45 & 45 & 45 \\
    Starch (corn starch)& 45 & - & 45 & 45 & 45 \\
    Cranberry (dried cranberry) & 10 & 10 & 20 & 10 & - \\
    Raisin & 10 & 10 & - & 10 & 20 \\
    Gin & 10 & 10 & 10 & - & 10 \\
    Water & - & - & - & 10 & - \\
    Cream cheese & 80 & 80 & 80 & 80 & 80 \\
    \hline
  \end{tabular}
  \caption{Composition of biscuits}  
  \label{table:biscuitcompo}
\end{table}

The biscuits cooking method and the weights of 
the ingredients were decided 
by one of the authors (M.Y.) of this work, who is a dietitian,
by reference to $19$ cream cheese cookies and biscuits recipes
(see Appendix~\ref{subsec:appendix:referencerecipes}).
Cream cheese (Yotsuba Milk Products Co., Ltd., Hokkaido, Japan) was stirred,
and sugar (Fuji Nihon Seito Corporation., Tokyo, Japan) was added and mixed.
Then, gin (THE BOTANIST, Scotland, England) was added and mixed,
and minced dried cranberry (Tomizawa Shouten, Tokyo, Japan)
and raisin (Tomizawa Shouten, Tokyo, Japan) were added and stirred.
Subsequently, mixed whole grain wheat flour (Tomizawa Shouten, Tokyo, Japan)
and corn starch (Tomizawa Shouten, Tokyo, Japan) were added and stirred,
and this dough was cooled at $\qty{4}{\degreeCelsius}$ for $\qty{60}{\minute}$.
The dough was flattened out to a thickness of $\qty{3}{\mm}$,
and was cut out to a diameter of $\qty{50}{\mm}$ using the cookie cutter.
Then, these were baked at $\qty{180}{\degreeCelsius}$ for $\qty{18}{\minute}$,
and were cooled down to room temperature.
Biscuits were packed and sealed in a vacuum package,
and were stored at $\qty{20}{\degreeCelsius}$
for
$\qty{24}{\hour}$ for the sensory evaluation. 

\subsection{Sensory evaluation}

Five biscuits (Control, NCS, NR, NG, and NCB) were placed on white dishes
and randomly tasted by blinded subjects.
The sensory evaluation was conducted by
$19$ untrained non-expert Japanese male and female students
between the ages of $21.3\pm2.6$ years old.
Perception of the intensity of
color (-3 [bright] to +3 [dark]),
texture (-3 [hard] to +3 [crispy]),
sweetness and sourness (-3 [weak] to +3 [strong])
of the biscuits
were examined.
The individual preferences for the
aroma, color, texture, taste, and overall judgment (-3 [dislike] to +3 [like]),
and palatability (-3 [bad] to +3 [good])
of the biscuits
were also examined.
Additionally, the preference ranking was evaluated using a rank test for significance \citep{newell1987expanded}.
{This sensory evaluation was approved by the Human Ethics Committee of Graduate School of Human Development and Environment, Kobe University (approval number: 653-2).}

Statistical analysis for results of sensory evaluation were performed using
Excel 2019 (Microsoft Japan Co., Ltd., Tokyo, Japan)
and
EZR software (Saitama Medical Center, Jichi Medical University, Saitama, Japan) \citep{kanda2013investigation},
which is a graphical user interface for
R (The R Foundation for Statistical Computing, Vienna, Austria, version 4.0.3).
The differences between the biscuits were compared using the Tukey's honest significant difference (HSD) test.
The ranking test was analyzed using the Newell and MacFarlane tables \citep{newell1987expanded}.
A \textit{p}-value $< 0.05$ was considered significant.


\section{Results}
\label{sec:results}

\subsection{Topological data analysis of the recipe data}

\paragraph{Persistence diagrams}

\begin{figure}[h]
  \centering
  \captionsetup[sub]{font=normalsize}
  \begin{subfigure}[b]{0.45\textwidth}    
    \caption{}
    \label{subfig:pd1}
    \centering
    \includegraphics[width=\textwidth]{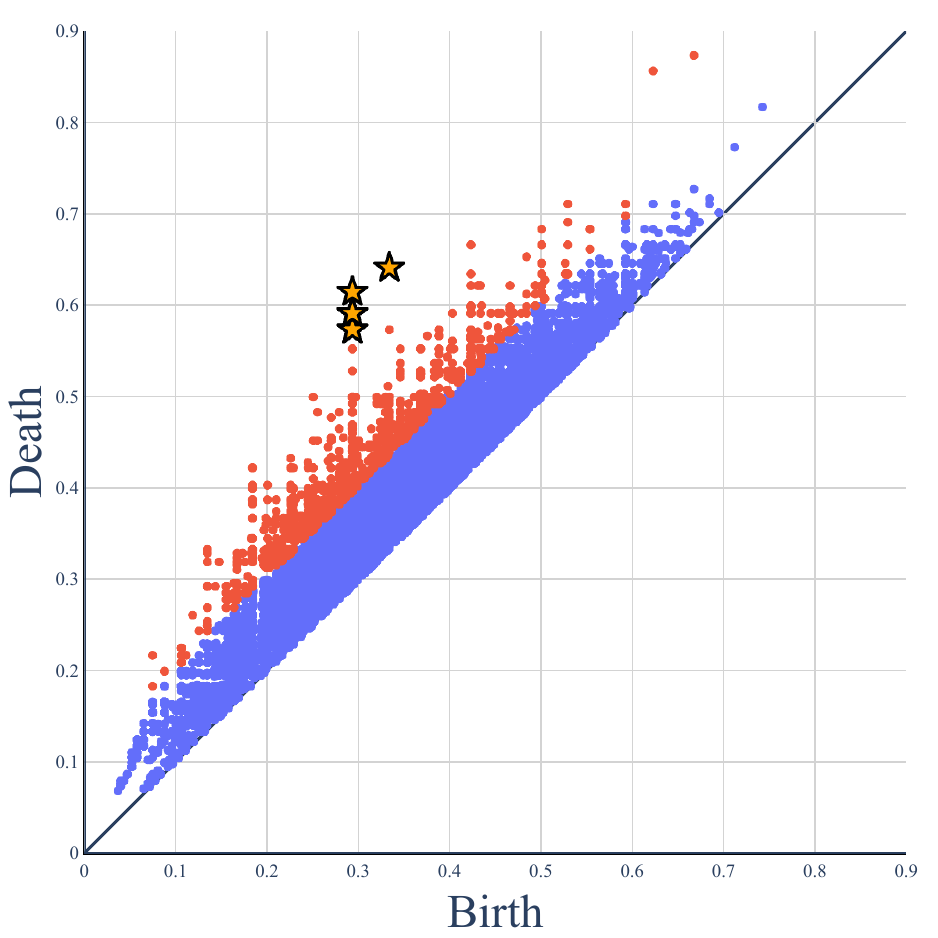}    
  \end{subfigure}
  \hfill
  \begin{subfigure}[b]{0.45\textwidth}
    \caption{}
    \label{subfig:pd1histo}
    \centering
    \includegraphics[width=\textwidth]{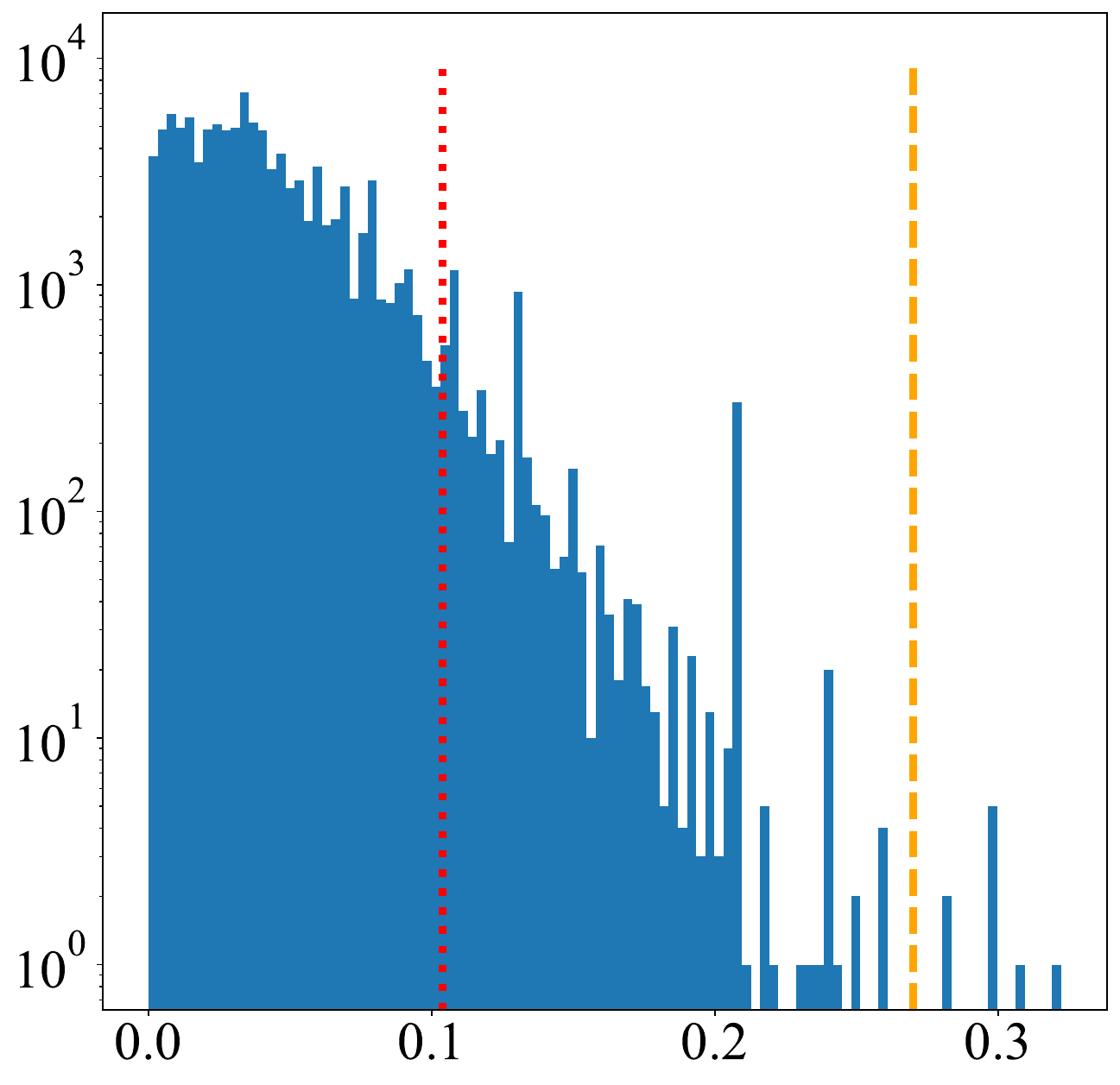}    
  \end{subfigure}
  \caption{
    Result of the persistent homology analysis
    of the recipe data from \citep{ahn2011flavor}.
    \subref{subfig:pd1},
    $1$st degree persistence diagram.
    The birth-death pairs with the top nine longest lifespans (lifespan $>0.27$) are starred.
    Out of all the birth-death pairs, $5\%$ (\datadep{$5,272$} birth-death pairs)
    with the longest lifespans (``the top $5\%$ birth-death pairs''), are shown in red.
    \subref{subfig:pd1histo},
    Histogram (with frequency in log-scale)
    of lifespans for the persistence diagram.
    The dissimilarity $0.27$ is marked with the dashed line.
    The dotted red line is at the smallest lifespan (\datadep{$\approx 0.104$}) of the top $5\%$ birth-death pairs.
  }
\end{figure}

We use the software ``Ripser'' \citep{Bauer2021Ripser} to compute persistent homology
and representative cycles\footnote{Using the branch ``representative-cycles'' of Ripser.}
for each birth-death pair in the persistence diagram.
The $1$st degree persistence diagram of the Vietoris-Rips filtration of
the recipe data with respect to the cosine dissimilarity
is shown in Figure~\ref{subfig:pd1}.
As suggested by the histogram of lifespans in Figure~\ref{subfig:pd1histo},
there are some birth-death pairs $(b,d)$ that
have longer lifespans compared to the other birth-death pairs (i.e. are far from the diagonal).
These are the birth-death pairs with lifespans greater than \datadep{$0.27$}, where a gap occurs in
the histogram (see also Table~\ref{table:bdl} in the Appendix).

\begin{figure}[h]
  \captionsetup[sub]{font=normalsize}
  \begin{subfigure}[b]{0.45\textwidth}
    \centering
    \caption{}
    \label{subfig:partiallistofrecipes}
    \begin{scriptsize}
      \begin{BVerbatim}
[almond,butter,cream,egg,vanilla,wheat]
[almond,butter,cream,cream_cheese,egg,vanilla,wheat]
[butter,cream,cream_cheese,egg,vanilla,wheat]
[butter,cream,cream_cheese,vanilla,wheat]
[butter,cream,vanilla,wheat]
[butter,cream,oat,vanilla,wheat]
[butter,cream,oat,wheat]
[butter,cream,wheat]
[butter,cream,gelatin,wheat]
[butter,cherry,cream,gelatin,wheat]
[cherry,cream,gelatin,wheat]
[cherry,cream,egg,gelatin,wheat]
[cherry,cream,egg,gelatin,orange,pineapple,wheat]
[cherry,cream,egg,gelatin,lemon,lime,milk,orange,pineapple,wheat]
[butter,cherry,cream,gelatin,lemon,lime,orange,pineapple,wheat]
[cherry,cream,gelatin,lemon,lime,orange,pineapple]
[cherry,lemon,lime,orange,pineapple]
[cherry,lemon,lime,orange,orange_juice,pineapple]
[cherry,lemon,lime,orange_juice,pineapple]
[lemon,lime,orange_juice,pineapple]
\end{BVerbatim}
 \end{scriptsize}
  \end{subfigure}
  \hfill
  \begin{subfigure}[b]{0.45\textwidth}    
    \centering
    \caption{}
    \label{subfig:regionalitytop}
    \includegraphics[width=0.87\textwidth]{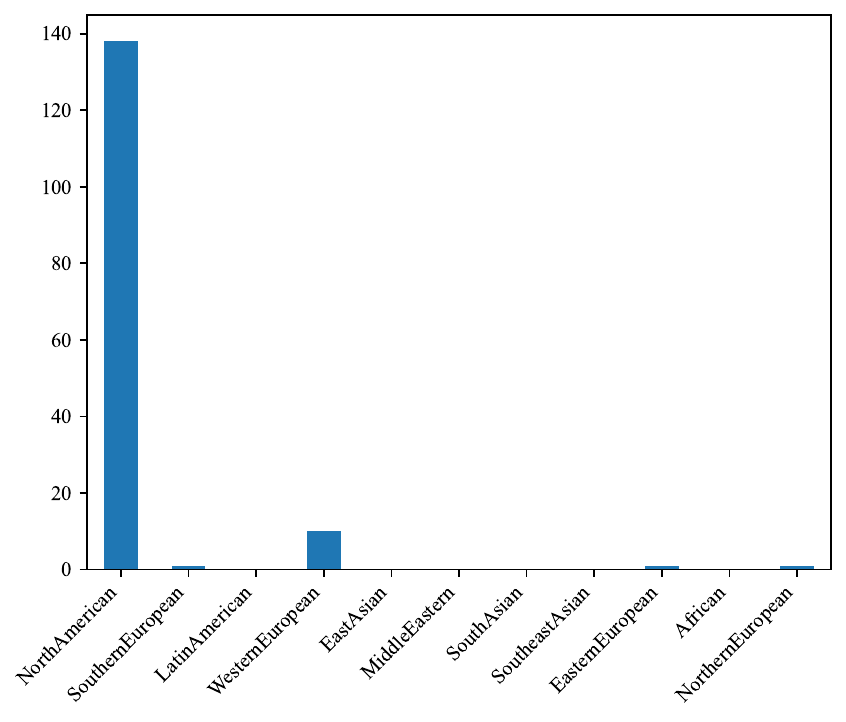}    
  \end{subfigure}
  \caption{Basic information about the representative cycle of the birth-death pair with the longest lifespan.
    \subref{subfig:partiallistofrecipes},
    partial list of the existing recipes (showing only $20$ out of a total of $97$ recipes)
    in a representative cycle for the birth-death pair with the longest lifespan.
    \subref{subfig:regionalitytop},
    regions associated to the $97$ recipes appearing in
    the representative cycle of the birth-death pair with the longest lifespan.
  }
  \label{fig:topbasic}
\end{figure}

As explained in the beginning of subsection~\ref{subsec:optimization},
each representative cycle $c$ is associated to a list of recipes $(r_1,r_2,\hdots, r_\ell)$
involved in its edges.
For example, in Figure~\ref{subfig:partiallistofrecipes},
we display a partial list of the recipes involved in
the cycle with the longest lifespan,
while List~\ref{list:cyclerecipes} in the Appendix are the recipes involved in
a cycle with the third longest lifespan.
We observed that (for the computed representative cycles with large lifespans)
going from one recipe $r_i$ to the next recipe $r_{i+1}$ in a cyclic fashion,
only one ingredient is added or removed;
that is, neighboring recipes in the detected cycles only differ by one ingredient.

In Figure~\ref{subfig:regionalitytop} we
display a bar plot of the regions associated to the \datadep{$97$} recipes appearing in
the representative cycle of the birth-death pair with the longest lifespan.
Note that due to the existence of multiple recipes with exactly the same
list of ingredients in the data, some ingredient combinations are associated to
multiple regions (possibly the same region multiple times), and thus the total number of regions exceeds
the number of recipes \datadep{$97$}.
We see that the representative cycle is mostly formed of recipes from North American cuisine
with some recipes from European cuisine mixed in.

In addition, we examine the top \datadep{$5,272$} birth-death pairs
($5\%$ of the total number, counting multiplicites)
with the longest lifespans (birth-death pairs colored red in Figure~\ref{subfig:pd1}).
Hereafter we shall use phrases such as ``top $5\%$ birth-death pairs'' or ``with the top $5\%$ longest lifespans''\footnote{Note that we need take into account multiplicity. The same birth-death pair may appear multiple times,
  and there may be ties in the lifespans.
  For example, simply ordering the \emph{distinct}
  lifespans as numbers
  (i.e.\ without repetition for the ties)
  and taking the top $5\%$ gives a different (larger) list of birth-death pairs.}
to refer to these birth-death pairs.

\begin{figure}[h!]
  \captionsetup[sub]{font=normalsize}
  \centering
  \begin{subfigure}[b]{0.45\textwidth}    
    \centering
    \caption{}
    \label{subfig:regionconcentration}
    \includegraphics[width=\textwidth]{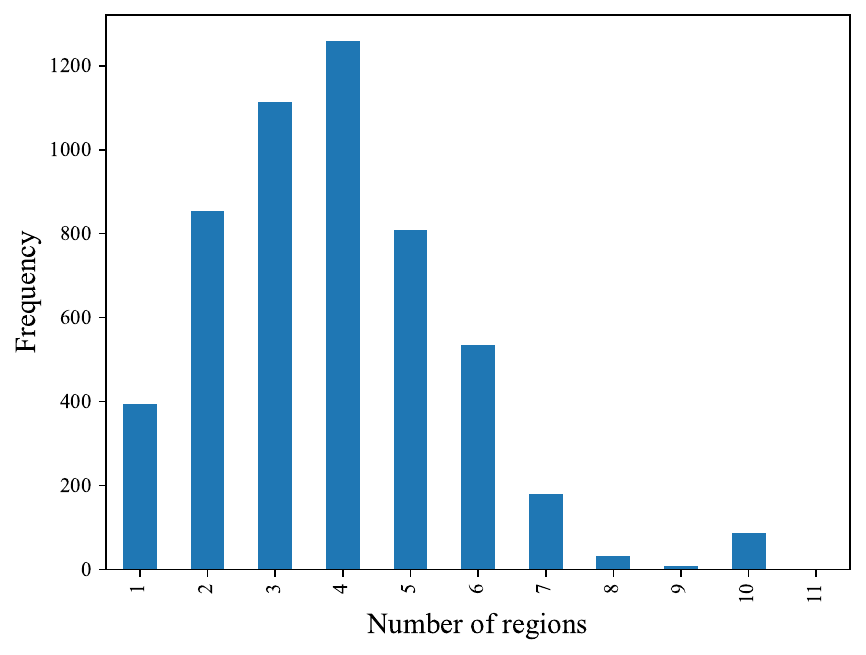}    
  \end{subfigure}
  \hfill
  \begin{subfigure}[b]{0.45\textwidth}    
    \centering
    \caption{}
    \label{subfig:regionalityaggregate}
    \includegraphics[width=\textwidth]{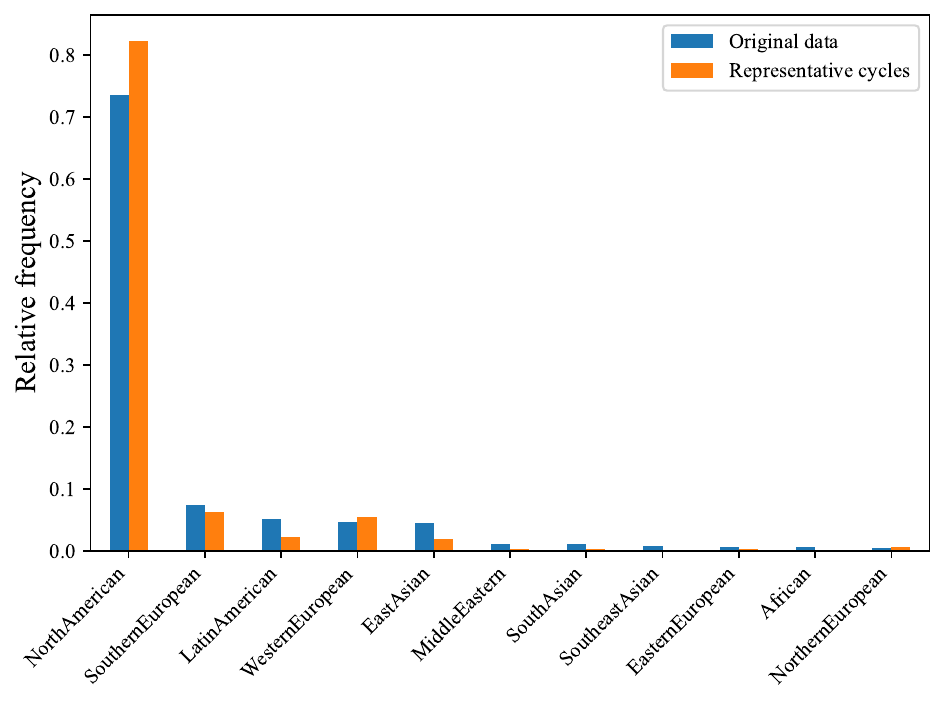}    
  \end{subfigure}
  \caption{
    Regionality of the
    representative cycles of the top $5\%$ birth-death pairs (\datadep{$5,272$} cycles).
    \subref{subfig:regionconcentration},
    histogram of number of distinct regions involved in the recipes of
    representative cycles of the top $5\%$ birth-death pairs.
    \subref{subfig:regionalityaggregate},
    relative frequencies of the regions associated to the recipes
    in the representative cycles
    of the top $5\%$ birth-death pairs, versus in those in the original data.}
  \label{fig:regionalityfivepercenter}
\end{figure}

We first examine the regionality of the recipes appearing
in the representative cycles of the top $5\%$ birth-death pairs.
In Figure~\ref{subfig:regionconcentration}, for these representative cycles,
we plot a histogram of the number of distinct regions involved, and see that most representative cycles involve recipes
from around \datadep{four} different regions.
On the other hand, plotting the relative frequencies of the regions of the representative cycles in aggregate,
versus the relative frequencies of the regions in the original data in
Figure~\ref{subfig:regionalityaggregate},
shows no significant difference (visually) in overall regions used in the representative cycles compared to the original recipe.

We then set up the combinatorial optimization problem described in
subsection~\ref{subsec:optimization} with $\nu=5$ in order to obtain,
for a birth death pair $(b,d)$ and its associated representative cycle $c$,
a set of solutions to the optimization problem~\eqref{equation:mastercombiopti}.
Each solution is a set of $\nu=5$ ingredients from the candidate set $S$
determined from $c$, such that the solution is as dissimilar as possible to existing recipes.
Note that there can be multiple solutions, representing ties in the dissimilarity to existing recipes.
For example, applied to the cycle with the longest lifespan,
we obtain solutions listed in List~\ref{list:somesolutionsfortopcycle}
(showing $5$ out of \datadep{635} total number of solutions).
\begin{mylist}
  \centering
  \begin{scriptsize}
\begin{BVerbatim}
(‘cranberry’, ‘cream cheese’, ‘gin’, ‘olive oil’, ‘raisin’)
(‘cranberry’, ‘cream cheese’, ‘gin’, ‘raisin’, ‘starch’)
(‘cranberry’, ‘cream cheese’, ‘gin’, ‘raisin’, ‘whole grain wheat flour’)
(‘cranberry’, ‘cream cheese’, ‘gin’, ‘starch’, ‘whole grain wheat flour’)
(‘cranberry’, ‘cream cheese’, ‘raisin’, ‘starch’, ‘whole grain wheat flour’)
\end{BVerbatim}
\end{scriptsize}
\caption{List of some solutions obtained to the combinatorial optimization problem associated to the
  representative cycles with longest lifespan.}
  \label{list:somesolutionsfortopcycle}
\end{mylist}
For the analysis below, we perform the computation for the representative cycles
of the top $5\%$ birth-death pairs,
but we only pick up a sample of up to $20$ solutions for each cycle,
as obtained by the applying the solver GLPK
(GNU Linear Programming Kit, Version 5.0,
\url{http://www.gnu.org/software/glpk/glpk.html})
to solve the optimization problem described in
subsection~\ref{subsec:optimization}.
Some cycles provide less than $5$ candidate ingredients; these provide no solutions.
Different cycles may lead to the same solutions; we remove duplicate solutions.
We obtain a total of \datadep{$31,478$} distinct solutions (suggested combinations) in total,
on which we perform further analysis below.

\paragraph{Combination novelty and ingredient usage statistics}

We first analyze the \datadep{$31,478$} suggestions from our method 
(consisting of up to $20$ solutions for the representative cycles
of the top $5\%$ birth-death pairs),
by checking whether or not they are present in the existing data as-is,
or as a sub-recipe of an existing recipe in the data.
We found that out of the \datadep{$31,478$} suggestions,
a small number \datadep{$61$ ($0.19\%$)} were in
fact existing combinations in the dataset,
and \datadep{$506$ ($1.6\%$)} were
strict subcombinations of existing recipes\footnote{Note that these two conditions are not mutually exclusive. A suggestion can be equal to an existing recipe, while at the same time being a strict subcombination of a (different) existing recipe.}.
For the former, this occurs when,
given the set of candidate ingredients $S$
computed from a representative cycle,
all combinations of specified size ($\nu = 5$ in our case) of $S$
are existing recipes.
Intuitively, this is more likely to happen
for a shorter representative cycle, which contains fewer recipes
and thus fewer total ingredients to use as candidates. In our collected sample of
suggestions, the largest lifespan where this occurs is \datadep{$\approx 0.150$}.
%
For the latter, we are checking whether or not
a solution $y_\ast$ (a suggested combination from our method)
satisfies the property that $y_\ast \subsetneq r$ for some existing recipe $r$ in the dataset.
In our collected sample of suggestions,
the largest lifespan where this occurs is 
\datadep{$\approx 0.207$}.
From the point of view of the novelty of the suggestions (with respect to the input data), the following interpretation can be given for these two possibilities.
Clearly, the former cases 
cannot be considered as novel since they are existing recipes as-is,
while the latter cases 
could be taken to mean removing one or more ingredients from an existing
recipe. Note that the majority  of the suggestions from our method
does not fall into either of these cases.

\begin{figure}[h]
  \centering  
  \includegraphics[width=0.8\textwidth]{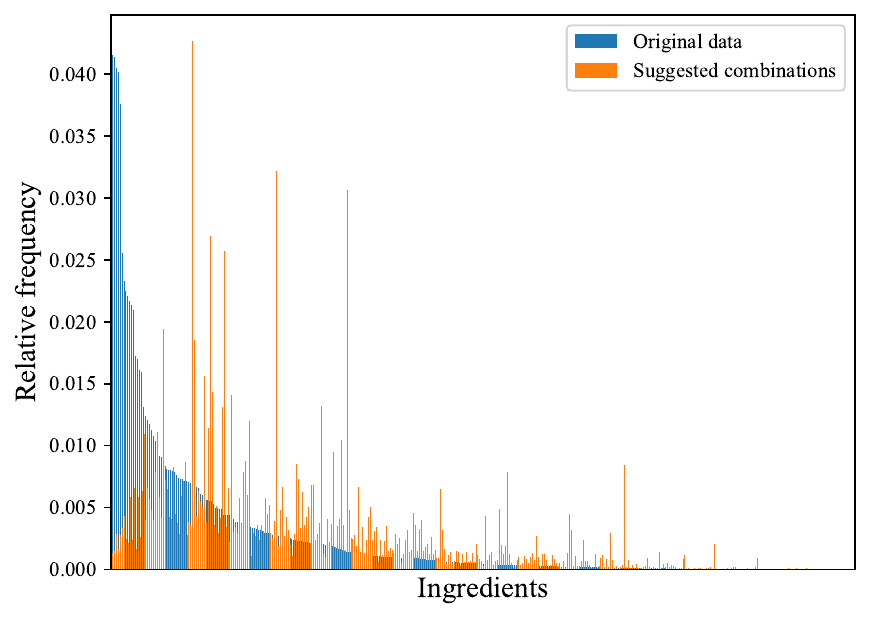}
  \caption{Relative frequencies of ingredients used in the original data (blue)
    versus in the suggestions obtained from our method (orange). We only apply our method to
    the representative cycles of the top $5\%$ birth-death pairs (\datadep{$5,272$} cycles),
    and only compute up to $20$ solutions for each (for a total of \datadep{$31,478$} solutions).
    Ingredients are ordered in decreasing frequency of usage in the original data.}
  \label{fig:ingredusecomparison}
\end{figure}

To further investigate the characteristics of the combinations of ingredients
suggested by our method, we explore ingredient usage statistics for both the original data and
the sample of \datadep{$31,478$} suggested combinations obtained from our method.
We plot the relative frequencies of ingredient usage in Figure~\ref{fig:ingredusecomparison}.
A visual inspection suggests that our method tends to suggest using ingredients not commonly used in the
original data.
This makes sense,
as we are maximizing dissimilarity of ingredient combinations compared to the
exisiting recipes. However, the method does not just preferentially use only rare ingredients
(the right tail of Figure~\ref{fig:ingredusecomparison}), because the initial analysis using
persistent homology constrains candidate ingredients to only those appearing in representative cycles.
Furthermore, it has been reported that
ingredient usage frequencies follow a power law distribution \citep{kinouchi_non-equilibrium_2008}.
In Appendix~\ref{subsec:appendix:powerlaw} we consider the distribution of the frequencies themselves
and performing curve fitting.

\subsection{Cooking and sensory evaluation for biscuits}
\label{subsec:cookingandsensory}

\begin{figure}[h]
    \centering
    \includegraphics[width=0.75\linewidth]{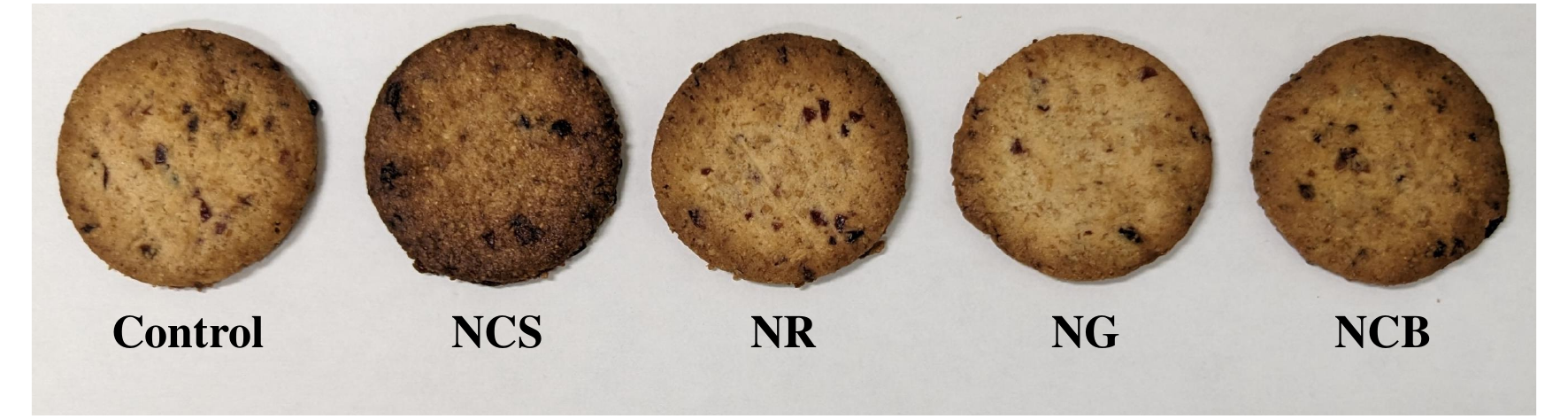}
    \caption{Visual appearance of biscuits }
    \label{fig:biscuits}
\end{figure}

We checked that suggested ingredient combinations from the above analysis
had viability as recipes of dishes.
For cooking,
we selected four similar ingredient combinations (List~\ref{list:somesolutionsfortopcycle})
obtained as solutions of our method.
These combinations (List~\ref{list:somesolutionsfortopcycle})
were not existing combinations in the dataset
nor strict subcombinations of existing recipes
\footnote{This can also be checked by noting the following. 
  The critical lifespan (considering lifespans in descending order)
  where our method started producing suggestions
  that were subcombinations of existing recipes is
  \datadep{$\approx 0.207$}.
  The lifespan of the cycle from which these suggestions were produced is \datadep{$> 0.27$},
  and thus cannot be existing combinations nor strict subcombinations.}.
%
%
%
Based on experience, the ingredient combinations were suggestive of ingredient lists for biscuits.
Therefore, we prepared biscuits using these lists, and evaluated their sensory properties.
The visual appearance of the resulting biscuits are displayed in Figure~\ref{fig:biscuits}.

\begin{figure}[h!]
  \captionsetup[sub]{font=normalsize}
  \centering
  \begin{subfigure}[b]{0.52\textwidth}    
    \centering
    \caption{}
    \label{subfig:sensoryA}
    \includegraphics[width=\textwidth]{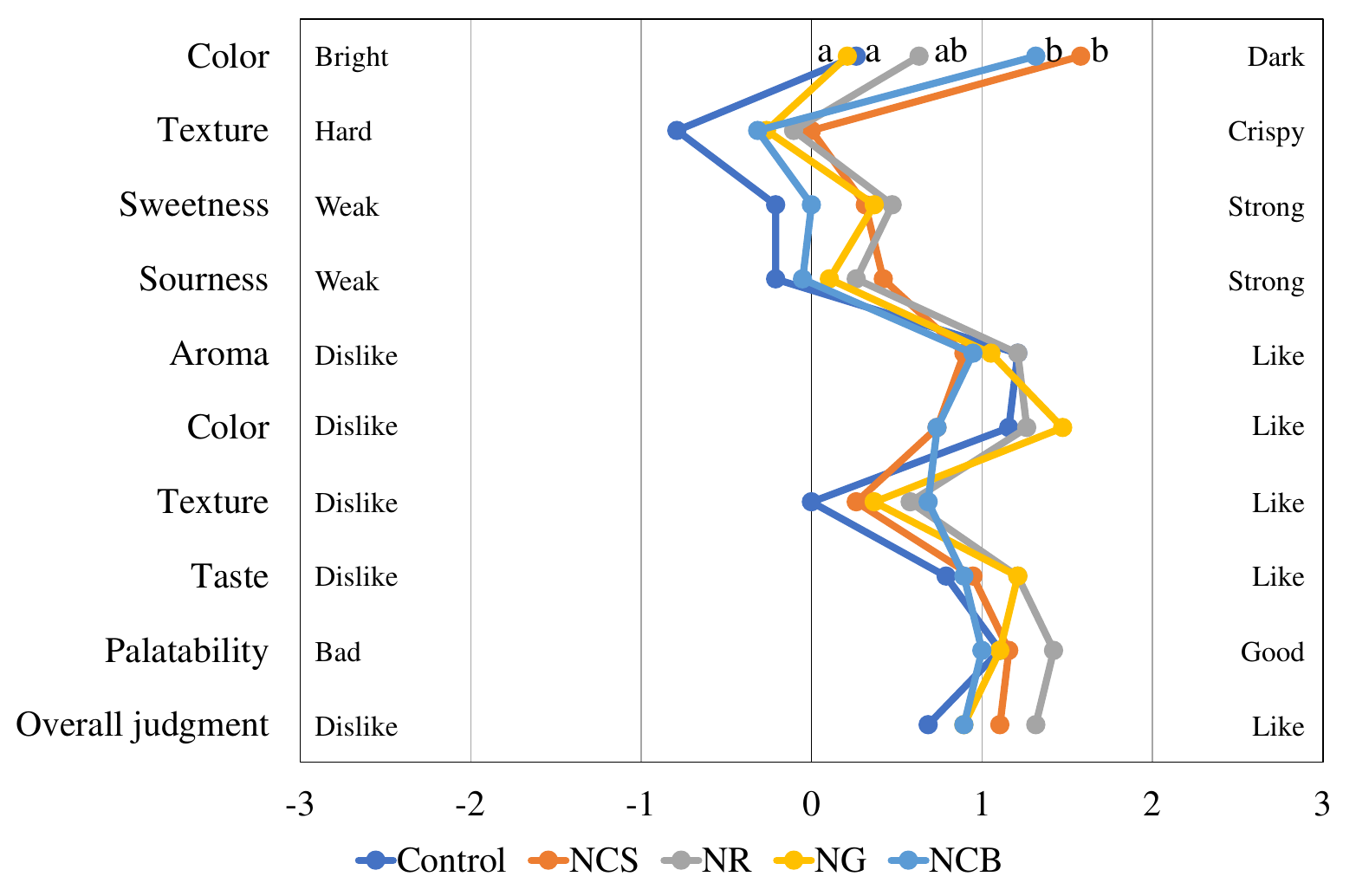}    
  \end{subfigure}
  \hfill
  \begin{subfigure}[b]{0.42\textwidth}    
    \centering
    \caption{}
    \label{subfig:sensoryB}
    \includegraphics[width=\textwidth]{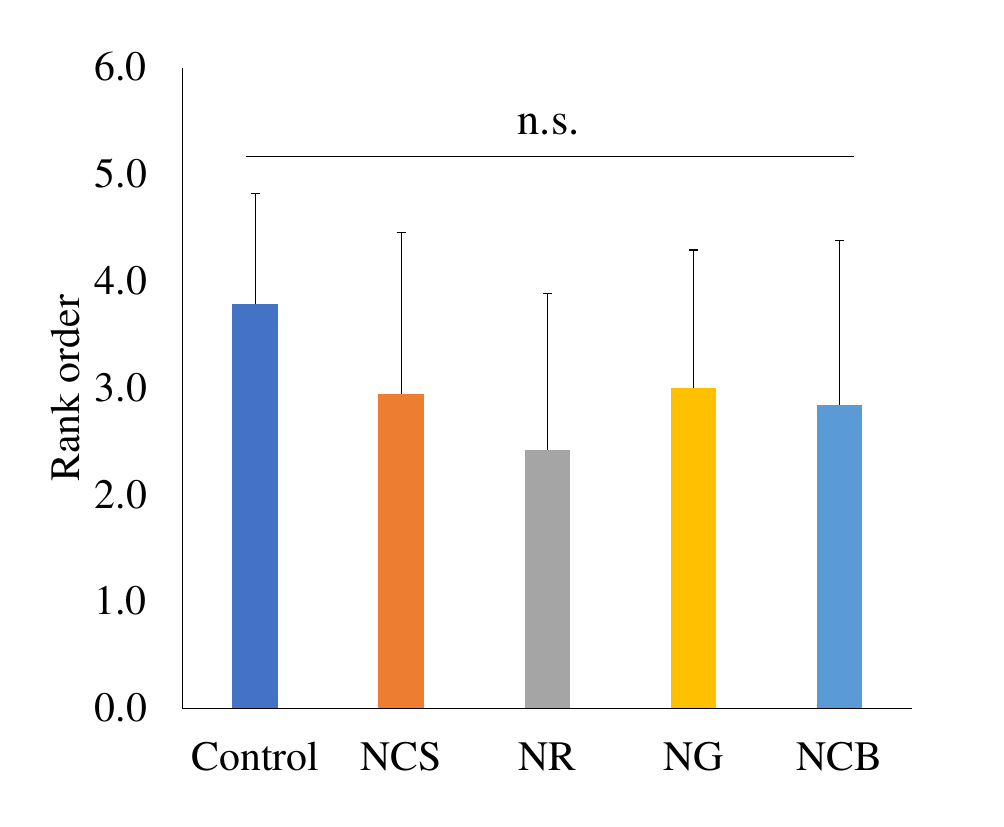}    
  \end{subfigure}
  \caption{
    The results of sensory evaluation. NCS, no corn starch; NR, no raisin; NG, no gin; NCB, no cranberry.
    \subref{subfig:sensoryA},
    sensory properties of biscuits, values are the mean ($n=19$). \textsuperscript{a-b} \textit{p }$< 0.05$ (Tukey’s HSD test).
    \subref{subfig:sensoryB},
    ranking order of biscuits analyzing the Newell and MacFarlane tables (Newell and MacFarlane, 1987), values are mean $\pm$ SD ($n=19$).}
  \label{fig:sensory}
\end{figure}

The results of the sensory evaluation are presented in Figure~\ref{subfig:sensoryA} and Figure~\ref{subfig:sensoryB}.
Statistical analysis shows that the color of NCS and NCB biscuits were judged to be darker than the other biscuits.
No statistically significant differences were found
for the other evaluations and the rank orders of the biscuits.
Scores of palatability and overall judgment of all biscuits were about $1.0$,
indicating that the suggested ingredient combinations
(List~\ref{list:somesolutionsfortopcycle} and Table~\ref{table:biscuitcompo})
are potentially viable for recipes of biscuits.

\section{Discussion}
\label{sec:discussion}

We discuss below some limitations of the current study.

In this work, for the data analysis we only considered recipes as
combinations of ingredients, without regard to cooking methods or ingredient weights, etc.
Indeed, we had to supplement the suggested combination from our method
with cooking methods and ingredient weights using references to existing recipes.
Furthermore, the method for vectorizing recipes uses a simple one-hot encoding of ingredients,
and thus similar ingredients (for example, cheddar cheese and cream cheese)
will be assigned different coordinates and
treated as dissimilar objects.

A priori, the following are some potential explanations
for the presence of holes is recipe space:
(1) missing data due to sampling issues; 
(2) ``recipes'' intentionally avoided,
due to the resulting dish being judged as not delicious, etc.; or
(3) novel recipes with culinary value.
The latter two possibilities are inherent features of the dataset.
Ideally one would like to be able to distinguish between the three cases.

Many solutions (ingredient combinations) were suggested by our method.
While in this work we cooked variations of cream cheese biscuits from some of the solutions,
further and more comprehensive investigation
into viability of the other combinations
for creating novel dishes is possible, especially
from the point of view of assisting chefs' creative cooking.
For example, the following combinations were also obtained from our method,
and could be interesting combinations to try to turn into novel dishes.
\begin{itemize}
\item cheddar cheese, melon, mustard, raisin, rum 
\item coffee, orange, pecan, shrimp, yeast
\item bell pepper, chicken liver, cocoa, cucumber, maple syrup
\end{itemize}

We are preparing follow-up work to address the above limitations.

In the sensory evaluation,
NCS and NCB biscuits were darker than other biscuits
(Figure~\ref{fig:sensory}A).
NCS biscuits did not include cornstarch, and
thus these were darker than other biscuits.
In addition, the sugar concentration of raisins is 
twice as much compared to dried cranberries \citep{swedishdatabase}.
NCB biscuits did not contain any dried cranberry and
contained twice the amount of raisin compared to other biscuits.
Therefore,
a potential explanation for why
NCB biscuits were darker than other biscuits
is a difference in the level of Maillard reaction.


\section{Conclusion}

We introduce the use of topological data analysis, especially persistent homology,
for the study of the space of culinary recipes.
In particular, we were able to identify sets of recipes surrounding the multiscale “holes”
in the space of existing recipes, which we exploit using combinatorial optimization
in order to generate novel ingredient combinations.
We perform analysis showing that the ingredient combinations suggested by our method are novel
with respect to the input data.
We also selected some ingredient combinations obtained from the analysis
and confirmed that we can cook variations of cream cheese biscuits from them.
A sensory evaluation study confirmed that these biscuits were acceptable enough.
Our findings indicate that topological data analysis has the potential
for providing new insights in the study of culinary recipes and
potentially for augmenting a chef’s creative process. 


\section*{Implications for gastronomy}

Recipes are important tools for any food and culinary culture, especially for preparing dishes in a delicious, nutritious, and safe manner. Chef’s culinary creativity is defined as the artistic and novel expression of the chef’s inner world, and the creativity is transferred to cuisine using their original recipes \citep{lee2021pate,lee2020creative}. However, creative chefs spend a lot of time creating novel recipes. In addition, the variety of dishes that a chef can create depends on the chef’s creativity, which can be subjective and potentially limited. A tool for suggesting novel ideas for dishes may be beneficial for augmenting a creative chef’s potential for novelty, and the use of tools from data science may contribute to this. In this study, we apply the topological data analysis in order to understand the “shape” of the space of cooking recipes, and “holes” are identified in the recipe space, and the potentiality of creating novel recipes is suggested. Our results will contribute to developing a tool for creating novel recipes and 
thus potentially contribute to the growth of gastronomy. 

\section*{Funding}
This research was supported by
the Casio Science Promotion Foundation (Grant number 40-29) and by the
Research Grant of Graduate School of Human Development and Environment, Kobe University.
E.G.E is supported by
  JSPS Grant-in-Aid for Transformative Research Areas (A) (22H05105) and 
  JSPS Grant-in-Aid for Scientific Research (C) (24K06846).

\section*{Acknowledgment}

None.

\section*{Author Statement}
\textbf{E.G.E.}:
Conceptualization,
Methodology,
Formal analysis,
Software,
Investigation - statistical and topological data analysis,
Writing,
Funding acquisition.
\textbf{Y.S.}:
Software,
Investigation - statistical and topological data analysis, Writing – review and editing.
\textbf{M.Y.}:
Conceptualization,
Methodology,
Investigation - cooking and sensory evaluation experiments,
Writing.

\section*{Declaration of competing interest}
The authors declare no conflict of interest.

\section*{Data availability}
Primary data is the recipe data from Supplementary Dataset 2 of \cite{ahn2011flavor},
available at \url{https://doi.org/10.1038/srep00196}.

\FloatBarrier

\bibliographystyle{elsarticle-harv} 
\bibliography{cas-refs}
 
 \newpage
 \appendix
\gdef\thesection{\Alph{section}}
\makeatletter
\renewcommand\@seccntformat[1]{\appendixname\ \csname the#1\endcsname.\hspace{0.5em}}
\makeatother

\section{Detailed definitions}
\label{sec:defns}

The details for the definition of homology are as follows.
First, for mathematical simplicity, we consider homology (and vector spaces) with coefficients in $\mathbb{F}_2$;
this is the field with two elements $\{0,1\}$ and with addition $1+1=0$.
Let $C_q(K)$ be the vector space of formal finite sums of $q$-simplices of $K$,
with coefficients in $\mathbb{F}_2$.
Letting the $q$-simplices of $K$ be $\{\sigma_1,\sigma_2,\hdots, \sigma_{\ell_q}\}$,
each element $c$ of $C_q(K)$ (called a \emph{$q$-chain}) is of the form
$
  c = \sum_{i=1}^{\ell_q} \delta_i \sigma_i
$
for some $\delta_i \in \mathbb{F}_2$.
Geometrically, $c$ can be thought of as the collection of $q$-simplices $\sigma_i$
for $i$ with $\delta_i \neq 0$ (i.e.\ $\delta_i = 1$).
The \emph{boundary} of a $q$-simplex $\sigma = \{v_0, v_1,\hdots, v_q\} \in K$ is
the $(q-1)$-chain
$
  \partial_q(\sigma) = \sum_{i=0}^q \{v_0,v_1,\hdots, \hat{v_i}, \hdots, v_q\}
$
where $\hat{v_i}$ means to remove the vertex $v_i$,
i.e.\ $\{v_0,v_1,\hdots, \hat{v_i}, \hdots, v_q\}$ is the $(q-1)$-simplex
with vertices $v_0,v_1,\hdots,v_q$, excluding $v_i$.
Defined this way, the above extends to a linear map $\partial_q$
from $C_q(K)$ to $C_{q-1}(K)$ called the \emph{$q$th boundary operator} of $K$.

The kernel
$
  Z_q(K) = \ker \partial_q = \{z \in C_q(K) \mid \partial_q(z) = 0\}
$
is called the space of $q$-cycles of $K$,
while the image
$
  B_q(K) = \im \partial_{q+1} = \{\partial_{q+1}(c) \mid c \in C_{q+1}(K)\}
$
is called the space of $q$-boundaries of $K$.
From the property that $\partial_q \partial_{q+1} = 0$,
it follows that $B_q(K) \subseteq Z_q(K)$. 
The quotient vector space
$
  H_q(K) = Z_q(K) / B_q(K)
$ is
called the \emph{$q$th homology group} of $K$.
For a $q$-cycle $z \in Z_q(K)$, its equivalence class $z + B_q(K)$ in $H_q(K)$
is also denoted by $[z] = z + B_q(K)$ and is called the \emph{homology class} of $z$.
On the other hand,
for a homology class $[z]$, $z$ is called a \emph{representative} of the homology class $[z]$.
Note that for a homology class $[z]$, any $z' \in [z]$ (i.e.\ $z' = z + b$ for some $b \in B_q(K)$)
can serve as a representative, since in this case $[z] = [z']$.

For persistent homology, the following
restatement of
the main theorem of \citep{edelsbrunner2002topological,carlsson2005computing}
(see also \citep[Section~2.1]{obayashi2018volume}, \citep[Theorem~2.6]{de2011dualities})
is one way of expressing the persistence diagram
and a set of its representative cycles, applied to Vietoris-Rips filtrations.
Below, for a $q$-cycle $c \in Z_q(V_t(X))$,
we let the homology class of $c$ in $H_q(V_t(X))$ be denoted by $[c]_t$ instead of just $[c]$,
to specify the threshold $t$ at which we are considering it.
\begin{theorem}
  \label{theorem:pd}
  There exists a set $q$-cycles $\{c_i\}_{i=1}^s$ and a unique multiset of pairs
  $\{(b_i,d_i)\}_{i=1}^s$ ($b_i \in \mathbb{R}$, $d_i \in \mathbb{R} \cup \{\infty\}$) such that the following hold.
  \begin{enumerate}
  \item If $t < b_i$ then $c_i \not\in Z_q(V_t(X))$.\label{ph:b1}
  \item If $b_i \leq t$ then $c_i \in Z_q(V_t(X))$.\label{ph:b2}
  \item If $b_i \leq t < d_i$ then $[c_i]_t \neq 0$.\label{ph:alive}
  \item If $d_i \leq t$ then $[c_i]_t = 0$.\label{ph:d}
  \item For each $t \in \mathbb{R}$,
    $\{[c_i]_t \mid i \text{ satisfies } b_i \leq t < d_i\}$ forms a basis for $H_q(V_t(X))$. \label{ph:basis}
  \end{enumerate}
\end{theorem}

\section{Additional analysis and results}
\label{sec:sample:appendix}

\subsection{Mean and standard deviation of the dissimilarities}
\label{subsec:appendix:dissim}

We further explore the cosine dissimilarities of the recipe data,
by considering the following random model.
A \emph{random bitstream} \citep{giller2012statistical} is a vector $B$ in $\{0,1\}^M$ where each component is independently
drawn from the Bernoulli distribution with parameter $p$.
We write this as $B \sim \bitstream{M}{p}$. Then, for $A \sim \bitstream{M}{p}$ and
$B \sim \bitstream{M}{q}$ independent,
the following approximations were given \citep{giller2012statistical}:
\begin{align}
  \mathbb{E}(\dcos(A,B)) & \approx 1 - \sqrt{pq} \label{eq:dcos_expectation}\\
  \mathrm{Var}(\dcos(A,B)) & \approx \frac{4-3p-3q + 2pq}{4M}. \label{eq:dcos_variance}
\end{align}

Using
$p = q = \frac{\text{Ave.~num.~ingreds.~per recipe}}{\text{Number of ingreds.}}
\approx \frac{\datadep{8.4936}}{\datadep{381}}
\approx \datadep{0.0223}$,
we obtained the numbers in Table~\ref{table:dcosstats}~column~(a).
The numbers in Table~\ref{table:dcosstats}~column~(b)
are computed using all pairwise dissimilarities of the data.
Since this violates the
independence assumption for $A$ and $B$ in $\dcos(A,B)$,
we also compute the mean and standard deviation of
the cosine dissimilarities between the recipes randomly paired without replacement
and give the result in
Table~\ref{table:dcosstats}~column~(c).
Compared to the (approximate) theoretical values of the random model,
the dissimilarities of recipes in the recipe data have smaller mean and larger standard deviation.
This suggests that the random bitstream model (which is purely random)
is not a good fit for the recipes.
Indeed, a copy-mutate model has been proposed to model cuisine evolution
\cite{kinouchi_non-equilibrium_2008}.
\begin{table}[h]
  \centering
  \begin{tabular}{r|c|c|c}
    & \thead[b]{(a)\\ from Eqs.~\eqref{eq:dcos_expectation}\eqref{eq:dcos_variance}}
    & \thead[b]{(b)\\ all pairwise\\ dissimilarities\\ of recipes}
    & \thead[b]{(c)\\ dissimilarities of\\ randomly paired  \\ recipes} \\\hline
            Mean $\dcos$  & \datadep{0.9777}  & \datadep{0.8681} & \datadep{0.8643}\\
    Std.~dev $\dcos$      & \datadep{0.05037} & \datadep{0.1414} & \datadep{0.1426}
  \end{tabular}
  \caption{Statistics for cosine dissimilarities}
  \label{table:dcosstats}
\end{table}

\FloatBarrier
\subsection{Birth-death pairs and representative cycles}
\label{subsec:appendix:cycles}
\begin{table}[h!]
  \centering
  \begin{tabular}{r|c|c|c|c}
    \thead[b]{~}
    & \thead[b]{birth-death pair}
    & \thead[b]{lifespan}
    & \thead[b]{num recipes in\\ rep.~cycle}
    & \thead[b]{num ingreds.~in\\ rep.~cycle} \\ \hline
    1 & (0.292893, 0.6151)   &  0.322207 &  97 & 40 \\
    2 & (0.333333, 0.641431) &  0.308098 &  64 & 27 \\
    3 & (0.292893, 0.591752) &  0.298859 &  85 & 41 \\
    4 & (0.292893, 0.591752) &  0.298859 &  35 & 20 \\
    5 & (0.292893, 0.591752) &  0.298859 &  83 & 34 \\
    6 & (0.292893, 0.591752) &  0.298859 &  77 & 34 \\
    7 & (0.292893, 0.591752) &  0.298859 & 125 & 53 \\
    8 & (0.292893, 0.573599) &  0.280706 &  99 & 41 \\
    9 & (0.292893, 0.573599) &  0.280706 &  84 & 36 \\
    10 & (0.292893, 0.552786) & 0.259893 &  93 & 40 \\
  \end{tabular}
  \caption{Birth-death pairs and number of recipes and ingredients in the representative cycles
    for the top 10 largest lifespans.}
  \label{table:bdl}
\end{table}

\begin{mylist}[h!]
  \centering
\begin{scriptsize}
\begin{BVerbatim}
['bean', 'beef', 'cayenne', 'cumin', 'garlic', 'onion', 'tomato']
['beef', 'cayenne', 'cumin', 'garlic', 'onion', 'tomato']
['beef', 'cayenne', 'garlic', 'onion', 'tomato']
['cayenne', 'garlic', 'onion', 'tomato']
['cayenne', 'garlic', 'olive_oil', 'onion', 'tomato']
['garlic', 'olive_oil', 'onion', 'tomato']
['olive_oil', 'onion', 'tomato']
['olive_oil', 'onion', 'pepper', 'tomato']
['olive_oil', 'onion', 'pepper']
['olive_oil', 'onion', 'pepper', 'pork']
['onion', 'pepper', 'pork']
['onion', 'pork']
['chicken', 'onion', 'pork']
['chicken', 'pork']
['pork']
['pork', 'sherry']
['sherry']
['cheese', 'sherry']
['cheese']
['bacon', 'cheese']
['bacon', 'cheese', 'pork']
['bacon', 'bean', 'cheese', 'pork', 'white_bread']
['bacon', 'bean', 'bread', 'cheese', 'onion', 'pork']
['bacon', 'bean', 'cheese', 'mushroom', 'onion']
['bean', 'mushroom', 'onion']
['mushroom', 'onion']
['beef', 'mushroom', 'onion']
['beef', 'mushroom', 'onion', 'tomato']
['beef', 'onion', 'tomato']
['beef', 'onion', 'tomato', 'vegetable']
['beef', 'black_pepper', 'onion', 'tomato', 'vegetable']
['beef', 'black_pepper', 'garlic', 'onion', 'red_kidney_bean', 'tomato', 'vegetable']
['beef', 'black_pepper', 'cayenne', 'garlic', 'onion', 'red_kidney_bean', 'tomato']
['bean', 'beef', 'black_pepper', 'cayenne', 'garlic', 'onion', 'tomato']
['bean', 'beef', 'black_pepper', 'cayenne', 'cumin', 'garlic', 'onion', 'tomato']
\end{BVerbatim}
\end{scriptsize}
\caption{List of recipes in an example of a representative cycle.}
\label{list:cyclerecipes}
\end{mylist}

\begin{figure}[h!]
  \centering
  \begin{subfigure}[b]{0.4\textwidth}
    \centering    
    \includegraphics[width=\textwidth]{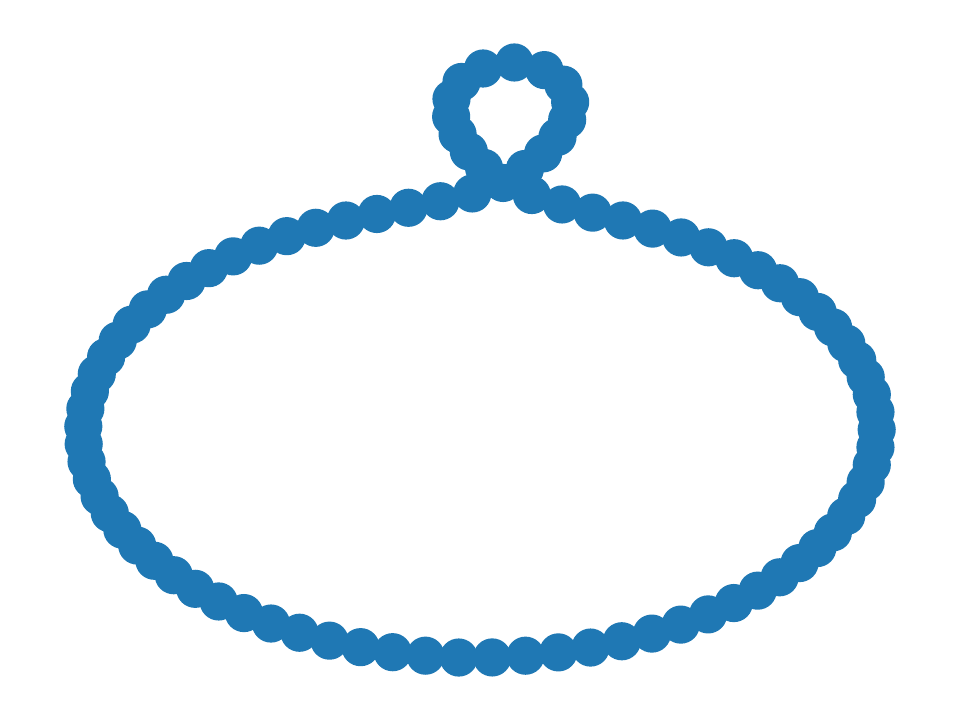}    
  \end{subfigure}
  \begin{subfigure}[b]{0.4\textwidth}
    \centering    
    \includegraphics[width=\textwidth]{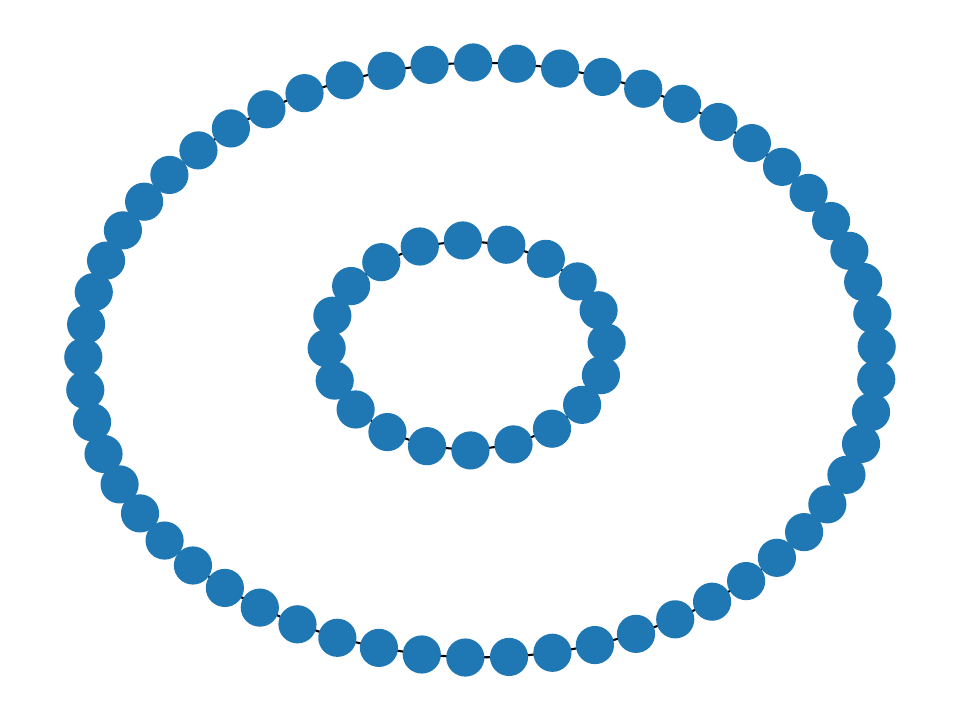}    
  \end{subfigure}
  \caption{Two examples of representative cycles that are not simple cycles.
  }
  \label{fig:nonsimple}
\end{figure}

\FloatBarrier
\subsection{Power law fitting}
\label{subsec:appendix:powerlaw}
\begin{figure}[h]
  \captionsetup[sub]{font=normalsize}
  \centering
  \begin{subfigure}[b]{0.45\textwidth}
    \centering
    \caption{}
    \includegraphics[width=\textwidth]{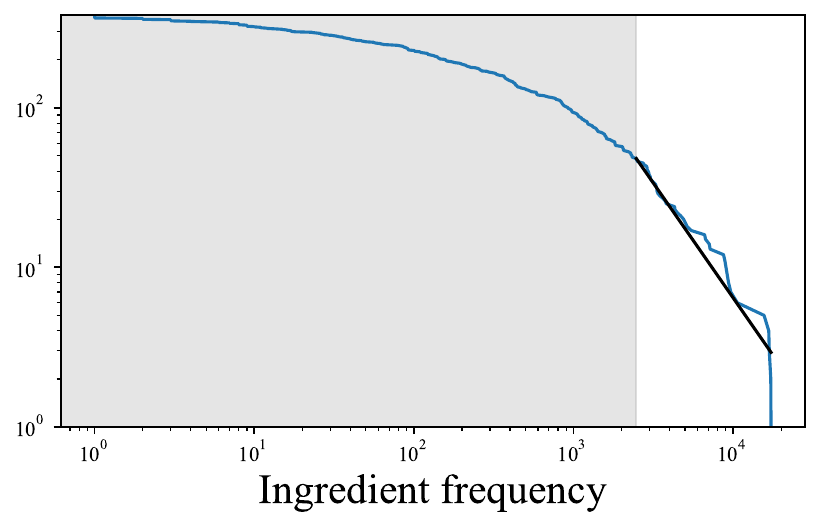}    
  \end{subfigure}
  \begin{subfigure}[b]{0.45\textwidth}
    \centering
    \caption{}
    \includegraphics[width=\textwidth]{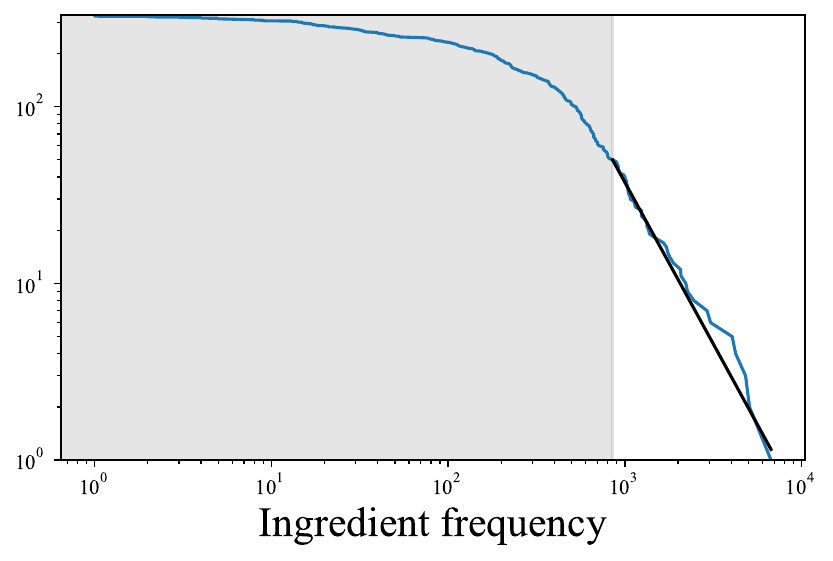}    
  \end{subfigure}
  \caption{Rank/frequency plots \citep{zipf1949human}
    in log-log scale for the ingredients used in the recipes of the
    (A)
    original data and
    (B)
    the combinations suggested by our method 
    consisting of up to $20$ solutions for the representative cycles with the top $5\%$ longest lifespans.
    Power law fitting computed using the python package ``powerlaw''~\citep{alstott_powerlaw_2014}.
    That is, we fit a curve $p(x) \sim x^{-\alpha}$ for $x \geq x_{\text{min}}$ where $x$ is frequency.
    The regions $x <  x_{\text{min}}$, which are excluded from the power law fit are shaded gray.
    A, 
    fitted parameters \datadep{$\alpha=2.438$, $x_{\text{min}} = 2,476$}.
    B, 
    fitted parameters \datadep{$\alpha=2.827$, $x_{\text{min}} = 852$}.
  }
  \label{fig:rankfreqfit}
\end{figure}

It has been reported that ingredient usage frequencies follow a power law distribution \citep{kinouchi_non-equilibrium_2008}. 
In such cases, the rank/frequency plot \citep{zipf1949human} appears as a straight line when plotted on log-log scale.
We perform power law fitting and show the results in Figure~\ref{fig:rankfreqfit}.
\FloatBarrier
\subsection{References for recipes }
\label{subsec:appendix:referencerecipes}

Biscuits cooking method and the weights of the ingredients
were decided by reference to the $19$
cream cheese cookie and biscuit recipes below.
References are all in Japanese,
and translations of recipe names and book titles
are unofficial and provided for reference only.

\newcommand*{\recipeinfo}[5]{
  \begin{CJK}{UTF8}{ipxm}#1\end{CJK}.
  \begin{CJK}{UTF8}{ipxm}#2\end{CJK} #3
  (#4)\\
  \href{#5}{#5}
  (Accessed 8 May 2022, in Japanese)
}

\begin{enumerate}[label={Recipe~\arabic*},wide=0pt,leftmargin=*]
\item
  \recipeinfo
  {Cookpad, 2014}
  {バターなし☆クリームチーズクッキー}
  {(recipe ID 2884806)}
  {Non-butter cream cheese cookies}
  {{https://cookpad.com/recipe/2884806}}

\item
  \recipeinfo
  {HomeBakery MARI no HEYA}
  {【レシピ】クリームチーズクッキー}
  {}
  {Recipe. Cream cheese cookies}
  {https://mari2.net/creamcheese-cookie/}  

\item
  \recipeinfo  
  {Nadia, 2022}
  {クリームチーズとくるみの米粉クッキー【バター不使用】}
  {(recipe ID 431995)}
  {Cream cheese and rice flour cookies with walnuts (non-butter)}
  {https://oceans-nadia.com/user/107837/recipe/431995} 
  
\item
  \recipeinfo
  {実験クッキング (Experimental cooking), 2016}
  {卵不使用米粉クリームチーズクッキー}
  {}
  {Non-egg rice flour and cream cheese cookies}
  {https://rdcooking.com/archives/393}  

\item
  \recipeinfo
  {DELISH KITCHEN}
  {しっとりさくさく！絞り出しクリームチーズクッキー}
  {}
  {Moisty and crispy! Squeezed cream cheese cookies}
  {https://delishkitchen.tv/recipes/171942334252974483}
  
\item
  \recipeinfo
  {手しごと私流 (Teshigotowatakushiryu), 2018}
  {ビニール袋で簡単！クリームチーズクッキー}
  {}
  {Easy cooking using a plastic bag! Cream cheese cookies}
  {https://watashiryuu.jugem.jp/?eid=176}
  
\item
  \recipeinfo
  {RAKUTEN recipe, 2011}
  {甘さ控えめ★クリームチーズクッキー レシピ・作り方}
  {(recipe ID 1650000839)}
  {Less sweetness. Cream cheese cookies recipe and instructions}
  {https://recipe.rakuten.co.jp/recipe/1650000839/}

\item
  \recipeinfo
  {RAKUTEN recipe, 2011}
  {卵いらず！簡単おやつ☆クリームチーズクッキー☆ レシピ・作り方}
  {(recipe ID 1480001331)}
  {Non-egg! Easy snack. Cream cheese cookies recipe and instructions}
  {https://recipe.rakuten.co.jp/recipe/1480001331/}

\item
  \recipeinfo
  {RAKUTEN recipe, 2013}
  {クリームチーズクッキー レシピ・作り方}
  {(recipe ID 1730008308)}
  {Cream cheese cookies recipe and instructions}
  {https://recipe.rakuten.co.jp/recipe/1730008308/}

\item
  \recipeinfo
  {RAKUTEN recipe, 2011}  
  {クリームチーズクッキー レシピ・作り方}
  {(recipe ID 1680000734)}
  {Cream cheese cookies recipe and instructions}
  {https://recipe.rakuten.co.jp/recipe/1680000734/}

\item
  \recipeinfo
  {RAKUTEN recipe, 2020}
  {簡単！クリームチーズクッキー♡ レシピ・作り方}
  {(recipe ID 1600033152)}
  {Easy! Cream cheese cookies recipe and instructions}
  {https://recipe.rakuten.co.jp/recipe/1600033152/}

\item
  \recipeinfo
  {RAKUTEN recipe, 2011}
  {超サックサク！簡単クリームチーズクッキー♪ レシピ・作り方}
  {(recipe ID 1130000693)}
  {Very crispy! Easy cream cheese cookies recipe and instructions}
  {https://recipe.rakuten.co.jp/recipe/1130000693/}

\item
  \recipeinfo
  {RAKUTEN recipe, 2012}
  {真っ赤なクランベリーのクリームチーズクッキー レシピ・作り方}
   {(recipe ID 1620006074)}
   {Red cranberry cream cheese cookies recipe and instructions}
   {https://recipe.rakuten.co.jp/recipe/1620006074/}

 \item
   \recipeinfo
   {RAKUTEN recipe, 2020}
   {クリームチーズクッキー レシピ・作り方}
   {(recipe ID 1230023044)}
   {Cream cheese cookies recipe and instructions}
   {https://recipe.rakuten.co.jp/recipe/1230023044/}

 \item
   \recipeinfo
   {RAKUTEN recipe, 2011}
   {洋酒の香りがおフランス　ガレット・ブルトンヌ レシピ・作り方}
   {(recipe ID 1140001988)}
   {French galette and bretonne with the aroma of Western liquor recipe and instructions}
   {https://recipe.rakuten.co.jp/recipe/1140001988/}

 \item S. Tamori, 2016.
   \begin{CJK}{UTF8}{ipxm}メープルビスケット\end{CJK}
   (Maple biscuits).
   In:
   \begin{CJK}{UTF8}{ipxm}米粉だから作れるとびきりおいしい焼き菓子\end{CJK}
   (Exceptionally delicious baked goods made with rice flour).
   IE-NO-HIKARI ASSOCIATION, pp. 48-51.
   (in Japanese)

 \item T. Okamura, 2020.
   \begin{CJK}{UTF8}{ipxm}かぼちゃの型抜きクッキー\end{CJK}
   (Squash cookies).
   In:
   \begin{CJK}{UTF8}{ipxm}卵, 牛乳, 白砂糖, 小麦粉なし. でも「ちゃんとおいしい」しあわせお菓子\end{CJK}
   (“Properly delicious” happy sweets that do not use eggs, milk, white sugar or flour).
   Kawade Shobo Shinsha., Ltd., pp. 52-53.
   (in Japanese)

 \item Y. Imai, 2020.
   \begin{CJK}{UTF8}{ipxm}2種のサブレ(紅茶レーズンとカモミールレモン)\end{CJK}
   (Two type of sables (tea \& raisin and chamomile \& lemon)).
   In:
   \begin{CJK}{UTF8}{ipxm}卵・乳製品・白砂糖をつかわない やさしいヴィーガン焼き菓子\end{CJK}
   (Easy vegan baked goods that do not use eggs, dairy products, or white sugar).
   Kawade Shobo Shinsha., Ltd., pp. 62-63.
   (in Japanese)
   
 \item T. Miyashita, 2009.
   \begin{CJK}{UTF8}{ipxm}レーズンドロップクッキー\end{CJK}
   (Raisin drop cookies).
   In:
   \begin{CJK}{UTF8}{ipxm}新調理学実習\end{CJK}
   (Text of new cooking practice).
   Dobunshoin Publishers, pp. 155.
   (in Japanese)
\end{enumerate}






\end{document}